\documentclass[10pt]{amsart}
\usepackage[T1]{fontenc}
\usepackage[utf8]{inputenc}

\usepackage[all]{xy}
\usepackage{amsmath}
\usepackage{amsthm}
\usepackage{caption}
\usepackage{graphicx}
\usepackage{multicol}
\usepackage{a4, amsfonts, amssymb, amsthm, latexsym}
\usepackage{enumerate}
\usepackage{fancyhdr}
\usepackage{color}
\usepackage{nicefrac}
\usepackage{multicol}
\usepackage{ textcomp }
\usepackage{ bbold }
\numberwithin{equation}{section}
\usepackage{mathrsfs,
     epsfig, multicol,anysize,graphicx,enumitem,mdwlist}

\newtheorem{thm}{Theorem}[section]
\newtheorem{lem}[thm]{Lemma}
\newtheorem{prop}[thm]{Proposition}
\newtheorem{cor}[thm]{Corollary}
\newtheorem{ex}[thm]{Example}
\newtheorem{defi}[thm]{Definition}

\newtheorem*{thm*}{Theorem}
\newtheorem*{prop*}{Proposition}

\def\im{\mathop{\mbox{\normalfont im}}}

\def\cok{\mathop{\mbox{\normalfont coker}}}

\DeclareMathOperator{\Hom}{Hom}
\DeclareMathOperator{\md}{mod}

\DeclareMathOperator{\Ext}{Ext}
\DeclareMathOperator{\Ind}{Ind}

\DeclareMathOperator{\End}{End}

\newtheorem*{defi*}{Definition}

\DeclareMathOperator{\cone}{cone}

\DeclareMathOperator{\Tr}{Tr}

\DeclareMathOperator{\Res}{Res}
\DeclareMathOperator{\colim}{colim}
\DeclareMathOperator{\U}{\bf U}
\DeclareMathOperator{\g}{\mathfrak{g}}
\DeclareMathOperator{\uu}{\bf u}

\newcommand{\mc}[1]{\mathcal{#1}} % Para las letras de las categorías.m_1\ell-2-m_0
\newcommand{\mb}[1]{\mathbb{#1}} %Para conjuntos numéricos
\newcommand{\mf}[1]{\mathfrak{#1}} %Para letraskazhdan - lusztig category góticas y de Lie
\newcommand{\fc}[3]{{#1}: {#2} \rightarrow {#3}}

\newcommand{\T}{\ensuremath{{\mathcal{T}}}}
\newcommand{\F}{\ensuremath{{\mathcal{F}}}}
\newcommand{\N}{\ensuremath{{\mathcal{N}}}}

\newcommand{\dem }{\noindent\textbf{Proof}. }
\newcommand{\findem }{\hfill $\Box$ \vskip0.5cm }
\marginsize{3cm}{3cm}{3cm}{3cm}

\def\im{\mathop{\mbox{\normalfont im}}}
\def\End{\mathop{\mbox{\normalfont End}}}

\def\cok{\mathop{\mbox{\normalfont coker}}}

\title{Derived counterparts of fusion categories of quantum groups}
\author{Juan Camilo Arias}
\date{ }
\address{Departamento de Matem\'{a}ticas, Universidad de los Andes, Carrera 1 N. 18A - 10, Bogot\'a, COLOMBIA. Tel: +571 3394999 ext. 3705 Fax: +571 3324427 }
\email{jc.arias147@uniandes.edu.co }
\thanks{The author was supported by {\it Convocatoria 2017-01 para la Financiaci\'on de Proyectos de Investigaci\'on Categor\'ia Estudiantes de Doctorado Candidatos. Proyecto: Grupos Cu\'anticos y Categor\'ias de Fusi\'on III }}
  \subjclass[2010]{Primary 17B37, 17B10, 18G99}
\keywords{Quantum groups, triangulated categories, fusion rings.}

\begin{document}

\begin{abstract}
In this paper, we study derived versions of the fusion category associated to Lusztig's quantum group $\U_q$. The categories that arise in this way are not semisimple but recovers the usual fusion ring when passing to Grothendieck rings. On the derived
level it turns out that it is possible to define fusion for $\U_q$ without using the notion of tilting modules. Hence, we arrive at a
definition of the fusion ring that makes sense in any spherical category. We apply this new definition to the small quantum group and we relate it with some rings appearing in \cite{AL}.
\end{abstract}

\maketitle
%\tableofcontents

\section{Introduction}

Let $\g$ be a semisimple finite dimensional complex Lie algebra and let $\U_q$ denote Lusztig's version of the quantized enveloping algebra at a root of unity $q$. Let $\mc{U}$ be the category of finite dimensional integrable $\U_q$-modules of type 1 and let $\mc{T} \subset \mc{U}$ be the full subcategory of tilting modules. The category $\T$ is self dual, closed under taking direct summands, finite direct sums and tensor products. We say that a tilting module is fusion if the highest weights of its indecomposable summands belong to the principal alcove of $\g$, in other case, if the highest weights of all its indecomposable summands are outside the principal alcove we say that the module is negligible. We denote by $\mc{N}$ and $\mc{F}$ the full subcategories of $\mc{T}$ whose objects are negligible and fusion tilting modules respectively. The category $\N$ is a tensor ideal of the category $\T$ and $\F$ is a additive subcategory of $\T$. As an additive category, $\T\cong \F\oplus \N$.\\

The fusion category of $\U_q$ is defined as the quotient category $\mc{T}/\mc{N}$, where $Obj(\T/\N)=Obj(\T)$ and for any $T_1, T_2\in \T$, $ \Hom_{\T/\mc{N}}(T_1, T_2) = \Hom_{\mc{T}}(T_1,T_2) /\mc{N}(T_1,T_2) $, where ${\mc{N}(T_1,T_2)} $ is the subgroup of morphisms form $T_1$ to $T_2$ that factors through $\N$. Since $\N$ is a tensor ideal of $\T$ the tensor product on $\T$ descends to a monoidal structure on $\T/\N$. As an additive category $\mc{T}/\mc{N}$ is canonically  equivalent to $\mc{F}$ but in order to see the tensor structure it must be regarded as a quotient category. The category $\mc{T}/\mc{N}$ is a rigid $\mb{C}$-linear tensor category. The fusion ring of the category $\mc{T}/\mc{N}$ is by definition the Grothendieck ring $\mc{R} = K_0(\mc{T}/\mc{N})$. It is a free abelian group of finite rank with basis in one to one correspondence with the irreducible fusion modules.\\

Denote by $K^b(Add)$ the bounded homotopy category of an additive category $Add$ and denote by $D^b(Ab)$ the bounded derived category of an abelian category $Ab$. Crucial to us is the important observation of Beilinson-Bezrukavnikov-Mirkovic, \cite{BBM}, that the canonical functor $\gamma: K^b(\mc{T}) \to D^b(\mc{U})$ is an equivalence of categories. Our first object of study is the Verdier quotient $K^b(\mc{T})/K^b(\mc{N})$. We notice that the canonical functor $\mc{T} \to \mc{T}/\mc{N}\cong\mc{F}$ induces a functor $K^b(\mc{T})/K^b(\mc{N}) \to K^b(\mc{F})$. It is not an equivalence of categories and $K^b(\mc{T})/K^b(\mc{N})$ is not a semisimple category. Nevertheless, as we show in Proposition \ref{Grothringquotient}, the Grothendieck rings of the categories $K^b(\mc{T})/K^b(\mc{N})$ and $\mc{T}/\mc{N}$ are isomorphic.\\

The equivalence $\gamma: K^b(\mc{T})\to D^b(\mc{U})$ induces an equivalence $K^b(\mc{T})/K^b(\mc{N}) \cong D^b(\mc{U})/\langle \mc{N} \rangle$, where $\langle \mc{N} \rangle \cong \gamma(K^b(\N))$  is the triangulated subcategory  of $D^b(\mc{U})$ generated by $\mc{N}$. We show that $\langle \mc{N} \rangle$ can be intrinsically described as the triangulated tensor ideal of $D^b(\mc{U})$ that is closed under retracts and generated by modules with singular highest weights, see Theorem \ref{NidealDU}. With this description, fusion rings can be defined without the notion of tilting modules, but it still depends on specifics of representation theory such us highest weight vectors.\\

We would like to give a description of the fusion category that only depends on the spherical structure. Let $\mc{N}_{\mc{U}} \subset \mc{U}$ be the category of \emph{all} negligible modules, that is, direct sums of indecomposable modules of quantum dimension $0$. We consider the category $D^b(\mc{U})/\langle \mc{N}_{\mc{U}} \rangle$, where $\langle \mc{N}_{\mc{U}} \rangle$ is the smallest triangulated subcategory containing $\mc{N}_{\mc{U}}$ and closed under retracts. In section \ref{conj}, we explain that the relation of the Grothendieck group $K_0(D^b(\mc{U})/\langle \mc{N}_{\mc{U}} \rangle)$ and the fusion ring $\mc{R}$ strongly depends if the retracts of the category $\langle \mc{N}_{\mc{U}} \rangle$ are negligible or not. We characterize this property in Proposition \ref{closedretractsequiv}. \\

Because the latter expression for the Grotendieck group makes sense in any spherical category, we suggest the following definition:

\begin{defi*}
 The fusion ring of an abelian spherical category $\mc{S}$ is $K_0(D^b(\mc{S})/\langle \mc{N}_{\mc{S}} \rangle)$ where $\mc{N}_{\mc{S}}$ is the full subcategory of negligible objects in $\mc{S}$, i.e., objects that are direct sums of indecomposable modules of quantum dimension zero and $\langle \mc{N}_{\mc{S}} \rangle$ is the smallest full triangulated subcategory that contains $\mc{N}_{\mc{S}}$ and is closed under retracts and tensor products with arbitrary modules. \end{defi*}

It would be interesting to calculate this ring for some spherical categories, for instance those arising from the spherical Hopf algebras of \cite{AAGITV}. In Proposition \ref{closedretractssmallsl2}, we show that for the case of the small quantum group $\uu_q$ of $\mf{sl}_2$ the retracts of the ideal $\langle \mc{N}_{\uu_q} \rangle$ are all modules belonging to $ \mc{N}_{\uu_q}$. This allow us to calculate in Proposition \ref{lachcomp}, using the above definition, the fusion ring of the small quantum group for the case of $\mf{sl}_2$ and show that it coincides with Lachowska's ring $\overline{Vr}$ of \cite{AL}. Here $Vr$ stands for Verlinde algebra which is another name for the fusion ring. Lachowska interprets $Vr$ as a quotient of the complexified Grothendieck ring of $\mc{U}$ and $\overline{Vr}$ is essentially obtained by basechanging $Vr$ to the complexified Grothendieck ring of the small quantum group. Thus, $D^b(\uu_q)/\langle \mc{N}_{\uu_q} \rangle$ is, at least for $\mf{sl}_2$, a categorification of $\overline{Vr}$. It would be very interesting to relate it with the other rings of Lachowska presented in \cite{AL}.\\

\noindent In addition, we study the stable category $\mathbf{S}(\mc{U}/\mc{N})$. In order to define it we use Beligiannis' theory of stabilization of left triangulated categories, \cite{B}. We first use an estimate based on the parabolic KL-polynomials (\cite{S}) to show that $\mc{N}$ is functorially finite, see Section \ref{contrfinitesubcat}. This implies that $\mc{U}/\mc{N}$ is a left triangulated category, which essentially means that there is an endofunctor (the so called shift functor) of the category $\mc{U}/\N$, which may not be invertible and a class of left distinguished triangles which satisfy similar axioms as the usual ones for a triangulated category. Formally inverting the shift functor we obtain the triangulated category $\mathbf{S}(\mc{U}/\mc{N})$. By Belangiannis' theory this category is equivalent to $K^{-,b}(\mc{N})/K^b(\mc{N})$ where $K^{-,b}(\mc{N})$ is the full subcategory of $K^{-}(\mc{N})$ whose objects are essentially $\mc{N}$-acyclic, hence, have bounded cohomology. Thus, the canonical functor $K^{-,b}(\mc{N}) \to D^b(\mc{U})$ induces a functor $$\mathbf{S}(\mc{U}/\mc{N}) \cong K^{-,b}(\mc{N})/K^b(\mc{N}) \to D^b(\mc{U})/K^b(\mc{N}).$$ We show that (unless $\g = \mathfrak{sl}_2$ in which case, all negligible modules are projective) this functor is not an equivalence of categories and the induced map on Grothendieck rings is surjective with a non-trivial kernel. Therefore, the ring $K_0(\mathbf{S}(\mc{U}/\mc{N}))$ can be thought of as an enhancement of the fusion ring $\mc{R}$. It would be interesting to explicitly calculate it in some cases.\\

\noindent This paper is organized as follows: Section 2 provides background material. In Section 3, we study the categories $K^b(\mc{T})/K^b(\mc{N}) \cong D^b(\mc{U})/\langle\mc{N}\rangle$ and $D^b(\mc{U})/\langle\mc{N}_{\mc{U}}\rangle$, and their Grothendieck rings. In Section 4, we study the stable category ${\bf S}(\mc{U}/\mc{N})$. Finally, in Section 5, we study fusion rings of more general spherical categories and compute the fusion ring of the small quantum group of $\mf{sl}_2$ and begin the task of relating our work to that of Lachowska \cite{AL}.

\section{Preliminaries}

In this section, we recall basic facts about quantum groups and tilting modules and we briefly review Belangiannis' theory of stabilization of left triangulated categories, \cite{B}.\\

\subsection{Root datum}
Let $\mf{g}$ be a finite dimensional semisimple complex Lie algebra. Let $\mf{h}\subseteq\mf{b}$ be a Cartan subalgebra contained in a Borel subalgebra of $\mf{g}$. Let $\Phi$ denote the corresponding root system and let $\Delta  = \{\alpha_i; \quad1 \leq i \leq n\}$ be the simple roots, so that the roots of $\mf{b}$ are positive. Let $Q \subset P \subset \mf{h}^*$ be the root lattice contained in the weight lattice. Let $P^{+}$ denote the dominant weights and let $Q^{+}$ be the semigroup generated by $\Delta$.  We equip $P$ with the partial order defined by $\mu \leq \lambda$ if and only if $\lambda - \mu \in Q^{+}$. For any root $\alpha\in \Phi$, we denote by $\alpha^{\vee} = \dfrac{2\alpha}{(\alpha,\alpha)}$ the corresponding coroot, where $(-,-)$ is the Killing form. We fix a non-negative integer $\ell$ which is prime to $3$ if $\g$ has components of type $G_2$.\\

The (finite) Weyl group is denoted by $W$. It is generated by the reflections $\fc{s_{\gamma}}{\mf{h}^{*}}{\mf{h}^{*}}$, $s_{\gamma}(\lambda) = \lambda - \langle \lambda,\gamma^{\vee} \rangle\gamma$,  through the hyperplanes $H_{\gamma} = \{ x\in \mb{R}\otimes_{\mb{Z}} Q \mid \langle x,\gamma^{\vee} \rangle =0 \}$, for $\gamma\in \Phi$, where $\langle \lambda,\gamma^{\vee} \rangle = 2(\lambda, \gamma)/(\gamma, \gamma)$. The longest element in $W$ is denoted by $\omega_0$.\\

For $\beta \in \Phi$ there exists $w\in W$ such that $\beta=w(\alpha_i)$ for some $i=1,\ldots, n$ (Theorem 10.3 \cite{H}). We set $\ell_{\beta} = \dfrac{\ell}{g.c.d(\ell, d_i)}$, where the $d_i$'s satisfy the relation $d_i\langle\alpha_i, \alpha_j^{\vee}\rangle = d_j\langle\alpha_j,\alpha_i^{\vee}\rangle$, for $i\neq j$ and $i,j=1,\ldots, n$. The {\it affine Weyl group} $W_{\ell}$, is the group generated by the reflections $\fc{s_{\beta,r}}{P}{P}$ through the affine hyperplanes $H_{\beta,r}=\{ x\in \mb{R}\otimes_{\mb{Z}}Q \mid \langle x+\rho, \beta^{\vee} \rangle = r \}$, $r\in\mb{Z}$,
defined for $\lambda\in P$ as $$ s_{\beta,r}\cdot\lambda = s_{\beta}\cdot\lambda + r\ell_{\beta}\beta  = \lambda - \langle \lambda + \rho, \beta^{\vee} \rangle\beta  + r\ell_{\beta}\beta$$ where $s_{\beta}\cdot \lambda = s_{\beta}(\lambda + \rho) - \rho$ is the dot action and $\rho$ is the half sum of the positive roots. $W_{\ell}$ is isomorphic to the semidirect product of the Weyl group $W$ and the translation group $\ell\mb{Z}\Delta$, i.e., $W_{\ell} \cong W \ltimes \ell\mb{Z}\Delta$.\\

We denote the principal alcove by $C_{\ell}=\{ \lambda\in P :0 < \langle \lambda + \rho , \alpha^{\vee} \rangle < \ell_{\alpha} \quad \forall \alpha\in \Phi^{+} \}$, and its closure by $\overline{C_{\ell}}=\{ \lambda\in P : 0 \leq\langle \lambda + \rho , \alpha^{\vee} \rangle \leq \ell_{\alpha} \quad\forall \alpha\in \Phi^{+} \}$.\\

\subsection{Quantized universal enveloping algebras}

We follow the notations of \cite{CP}. Let $v$ be an indeterminate, $\mc{A}:=\mb{Z}[v,v^{-1}]$ the Laurent polynomials with coefficients in $\mb{Z}$ and $\mb{Q}(v)$ its quotient field. The quantized universal enveloping algebra $U_v$ of the Lie algebra $\mf{g}$ is the associative algebra over $\mb{Q}(v)$ with generators $E_{i} = E_{\alpha_i}, F_i = F_{\alpha_i}$ and $K_i = K_{\alpha_i}$, for $ \alpha_i\in \Delta$, subject to the
relations:

$$K_{i}K_{j} = K_{j}K_i, \quad\quad K_{i}K_{-i} = K_0=1$$
$$K_{j}E_{i}K_{-j} = v^{d_i\langle \alpha_j,\alpha_i^{\vee} \rangle}E_{i}, \quad\quad K_{j}F_{i}K_{-j} = v^{-d_i\langle \alpha_j,\alpha_i^{\vee} \rangle}F_{i}$$
$$E_{i}F_{j} - F_{j}E_{i} = \delta_{ij}\frac{K_{i} - K_{-i}}{v^{d_i} - v^{-d_i}}$$
and certain quantum Serre relations that we do not recall here ($\delta_{ij}$ stands for the Kronecker delta). The algebra $U_v$ has a Hopf algebra structure, see \cite{CP} for details.\\

Lusztig's integral form $\U_{\mc{A}}$ is the $\mc{A}$-subalgebra of $U_v$ generated by the divided powers $E_i^{(N)} = E_i^N/[N]!$, $F_i^{(N)}=F_i^N/[N]!$, $K_i$ and $K_i^{-1}$ for $1\leq i\leq n$ and $N\geq 0$. Thus, multiplication gives an isomorphism $\U_{\mc{A}}\otimes_{\mc{A}} \mb{Q}(v)\cong U_v$.\\

Let $\U_{\mc{A}}^0$ the $\mc{A}$-subalgebra generated by $K_i$, $K_i^{-1}$ and the symbols $$ \begin{bmatrix} K_i; c \\ r \end{bmatrix}_{v_i} = \prod_{s=1}^{r} \frac{K_iv_i^{c+1-s} - K_i^{-1}v_i^{s-1-c}}{v_i^s-v_i^{-s}} $$ for all $i=1,2,\ldots, n$, $c\in \mb{Z}$ and $r\in\mb{N}$.\\

We fix from now on a primitive $\ell^{th}$-root of unity $q\in\mb{C}$. We consider the field of complex numbers $\mb{C}$ as an $\mc{A}$-module using the homomorphism $\mc{A}\rightarrow \mb{C}$, $v\mapsto q$. We define Lusztig's quantum group at a root of unity $q$ as,
$$\U_q  := \U_{\mc{A}}\otimes_{\mc{A}}\mb{C}.$$

The algebra $\U_q$ inherits a Hopf algebra structure from $U_v$. We define Hopf subalgebra $\U_q^{0}$ by $\U_{\mc{A}}^{0}\otimes_{\mc{A}}\mb{C}$.

\subsection{Categories of integrable $\U_q$-modules}

Let $M$ be an $\U_{\mc{A}}$-module. For a complex valued character $\lambda$ on $\U_{q}^{0}$ let $M^\lambda$ denote the corresponding weight space. $M$ is called integrable of type 1, if $M$ is the direct sum of its weight spaces and for all $x \in M$ there exists $r_x>0$ such that $E_i^{(r)}x = F_i^{(r)}x = 0, 1\leq i\leq n, r \geq r_x$. The latter condition automatically holds if $M$ is finite dimensional. We will denote the category of all integrable type 1 $\U_{q}$-modules by $\mc{U}^{int}$.\\

Let $\mc{U} \subset \mc{U}^{int}$ be the full subcategory of finite dimensional modules. For $M \in \mc{U}$ define the dual $M^{\ast}=\Hom_{\mb{C}}(M, \mb{C})$ with the action $(uf)(m)=f(\omega(S(u))x)$, for $f\in \Hom_{\mb{C}}(M, \mb{C})$, $u\in\U_q$ and $m\in M$. Here $S$ denotes the antipode of $\U_q$ and $\fc{\omega}{\U_q}{\U_q}$ is the Cartan involution, see Lemma 4.6 of \cite{J}. Note that $M$ and $M^*$ have the same formal character.\\

Let $\Delta(\lambda)$ denote the \emph{standard or Weyl module} of highest weight $\lambda$ and let $\nabla(\lambda) =
\Delta(\lambda)^*$, the \emph{costandard module}. Let $L(\lambda)$ denote the unique \emph{irreducible} quotient of
$\Delta(\lambda)$.  The category $\mc{U}$ has enough projectives and enough injectives. Moreover injective and projective modules coincide.
Let $I(\lambda)$ (resp. $P(\lambda)$) denote the \emph{injective hull} (resp. \emph{projective cover}) of $L(\lambda)$.
Note that $L(\lambda) \cong L(\lambda)^*$ and $P(\lambda) \cong P(\lambda)^* \cong I(\lambda) \cong I(\lambda)^*$. There is the {\it Steinberg module} $St := \Delta((\ell-1)\rho)$. It is irreducible, self dual and projective.\\

Two weights $\mu, \lambda \in P$ are said to be \emph{linked} if $\mu \in W_{\ell}\cdot \lambda$.  Let $\mc{U}^\lambda$ be the subcategory of $\mc{U}$ whose objects have composition factors $L(\mu)$ for $\mu$ linked to $\lambda$. By the linkage principle (Theorem 4.3 and Corollary 4.4 in \cite{A1}) we have the orthogonal decomposition $\mc{U} = \oplus_{\lambda \in P/W_{\ell}}\mc{U}^\lambda$.  In particular, any indecomposable module belongs to $\mc{U}^\lambda$ for some $\lambda$.  We have $\Delta(\lambda), \nabla(\lambda), L(\lambda), P(\lambda), I(\lambda) \in \mc{U}^\lambda$.

\subsection{Tilting modules}

 Let $M$ be a finite dimensional $\U_q$-module. $M$ has a \emph{standard} (resp. \emph{costandard}) \emph{filtration} if there exists a chain of
submodules $0  = V_0 \subset V_1 \subset \cdots \subset V_{p-1} \subset V_p = M $ such that $V_r / V_{r-1} \cong \Delta(\lambda_r)$ (resp. $V_r / V_{r-1} \cong \nabla(\lambda_r)$) for some $\lambda_r \in P^{+}$ and $r = 1, \ldots, p$. It is well known that $\Ext^{i}_{\mc{U}}(\Delta(\lambda), \nabla(\mu)) = 0$, for all $\lambda, \mu$ and $i>0$, \cite{HA}. From this it follows that if $M$ admits a standard filtration and $M'$ admits a costandard filtration then, $\Ext^{i}_{\mc{U}}(M, M') = 0$ for $i>0$.\\

\begin{defi} A finite dimensional $\U_q$-module is called \emph{tilting} if it has a standard filtration and a costandard filtration.\end{defi}

A standard module is tilting if and only if it is irreducible. In particular, standard modules with highest weight in the principal alcove are tilting.\\

Let $\mc{T}$ denote the full subcategory of $\mc{U}$ whose objects the tilting modules. By the above, all higher extensions between tilting
modules (calculated in $\mc{U}$) vanishes. Tilting modules are self-dual and they are closed under taking direct summands, finite direct sums and tensor products. The last fact is rather deep, see \cite{HA} and references therein. By Theorem 2.5 of \cite{HA}, for any dominant
weight $\lambda$ there exists a unique up to isomorphisms indecomposable tilting module $T(\lambda) \in \mc{U}^\lambda$. All
tilting modules are isomorphic to direct sums of such in a unique way up to permutations of factors. Thus, by the linkage principle there is a block decomposition $\mc{T} = \oplus_{\lambda \in P/W_{\ell}} \mc{T}^\lambda$ where $\mc{T}^\lambda = \mc{U}^\lambda \cap \mc{T}$.\\

For $\lambda\in P^{+}$, write $\lambda = \lambda_0 + \ell\lambda_1$, where $0\leq\langle \lambda_0, \alpha_i^{\vee} \rangle < \ell $ for all simple roots $\alpha_i$. Put $\overline{\lambda} := 2(\ell-1)\rho + w_0\lambda_0 + \ell\lambda_1$. Then $I(\lambda) = T(\overline{\lambda})$. In particular, any injective module is tilting. It is also known that any injective module is isomorphic to a direct summand in $St \otimes T$ for certain tilting module $T$, see Theorem 9.12 in \cite{APW}.

\subsection{The Fusion category}\label{fusioncat}

A fusion category is a rigid semisimple $\mb{C}$-linear monoidal category with finitely many isomorphism classes of simple objects, such that the unit object is indecomposable. We recall here the construction of the fusion category associated to the quantum group $\U_q$, \cite{AP}.\\

For $M\in \mc{U}$ and $f \in \End_{\mc{U}}(M)$, let $Tr_q(f) :=  \Tr(K_{2\rho}f)$ its \emph{quantum trace}, where
$K_{2\rho} = \prod_{\beta\in \Phi^{+}}K_{\beta}$ and $Tr$ is the usual trace of a $\mb{C}$-linear endomorphism.  The \emph{quantum dimension} of $M$ is $dim_q(M) := \Tr_q(1_M)$.\\

A module $M \in \mc{U}$ is {\it negligible} if $\Tr_q(f)=0$ for any $f \in \End_{\mc{U}}(M)$. An indecomposable module is negligible if and only if its quantum dimension is $0$. Hence, $M$ is negligible if and only if the quantum dimension of all its indecomposable direct summands is $0$. We denote by $\mc{N}_{\mc{U}} \subset \mc{U}$ the full subcategory of negligible modules. The category $\mc{N}_{\mc{U}}$ is a tensor ideal in $\mc{U}$.  It is  known that  $\Delta(\lambda) \in \mc{N}_{\mc{U}}$ if and only if $\lambda$ is a $\ell$-\emph{singular weight}, equivalently $\langle \lambda + \rho , \beta^{\vee} \rangle$ is divisible by $\ell$ for some positive root $\beta$.\\

Let $\mc{N} = \mc{N}_{\mc{U}}  \cap \mc{T} \subset \mc{T}$ denote the full subcategory of negligible tilting modules. The category $\mc{N}$ is a tensor ideal in $\mc{T}$. It is known that $T(\lambda) \in \mc{N}$ if and only if $\lambda \notin C_{\ell}$. All injective modules are negligible tilting modules. Unless $\g = \mathfrak{sl}_2$ there are negligible tilting modules which are not injective.\\

For $\lambda \in C_\ell$ we have $T(\lambda) = \Delta(\lambda) = L(\lambda)$. We refer to such a $T(\lambda)$ as an irreducible fusion module.
A fusion module is a module isomorphic to a direct sum of such. Let $\mc{F}$ be the full subcategory of $\mc{T}$ whose objects are fusion modules.
Thus, $\mc{F}$ is a semi-simple abelian category. A non-zero map between fusion modules cannot factor through a negligible tilting module. \\

Any  $T \in \mc{T}$ is non-canonically isomorphic to a direct sum $F \oplus N$, for $F\in\mc{F}$ and $N\in \mc{N}$. However, there is a way to construct the fusion part of $T$ canonically. It goes as follows, see \cite{CP} Proposition 11.3.18: Let $T^{\vee}\subset T$ be the maximal submodule of $T$ belonging to $\mc{F}$ and  let $T^{\wedge}$ be the maximal quotient of $T$ belonging to $\mc{F}$. Denote by $T^{\mc{F}}$ the image of $T^\vee$ under the projection $T \to T^\wedge$. Then $T^{\mc{F}}$ is isomorphic to the fusion part of $T$ and the assignment $T
\mapsto T^{\mc{F}}$ defines a functor \begin{equation}\label{TtoFpart} (\ )^{\mc{F}}: \mc{T} \to \mc{F}\end{equation} \noindent whose kernel is $\mc{N}$.\\

The category $\mc{F}$ is not closed under the tensor products of representations. Hence, since the tensor product of two representations of $\mc{F}$ is a tilting module, we equip $\F$ with the monoidal structure given by the reduced tensor product $$ F_1\overline{\otimes} F_2 := (F_1\otimes F_2)^{\mc{F}}. $$ The category $\F$ with $\overline{\otimes}$ is a fusion category and the fusion ring $\mc{R}$ is by definition the Grothendieck ring $K_0(\F)$.\\

The fusion category can also be thought of as a quotient. Consider the quotient category $\T/\mc{N}$ (see the introduction for its definition). Since $\N$ is a tensor ideal of $\T$ the tensor product on $\T$ descends to a monoidal structure (again called tensor product) on $\T/\N$. Note that the composition $\F \to \T \to \T/\mc{N}$ is an equivalence of additive categories which induces a monoidal structure on the category $\F$. This monoidal structure is equivalent to the one given by $\overline{\otimes}$.

\subsection{Verdier quotients and Grothendieck groups of triangulated categories}\label{GTI}  Let $\mc{C}$ and $\mc{D}$ be triangulated categories and $\fc{F}{\mc{C}}{\mc{D}}$ a triangulated functor. The kernel of $F$ is the full subcategory $\ker F$ of $\mc{C}$ whose objects maps to objects in $\mc{D}$ isomorphic to zero. Note that $\ker F$ contains all the direct summands of its objects. In general, we say that a subcategory $\mc{E}$ of $\mc{C}$ is thick or closed under retracts, if it is triangulated and it contains all direct summands of its objects.\\

Given a triangulated category $\mc{D}$ and a triangulated subcategory $\mc{C}$ (not necessarily thick) there exists a triangulated category $\mc{D}/\mc{C}$, and a triangulated functor $\fc{F_{univ}}{\mc{D}}{\mc{D}/\mc{C}}$ so that $\mc{C}$ is in the kernel of $F_{univ}$, and $F_{univ}$ is universal with the property: If $\mc{Q}$ is a triangulated category and $\fc{F}{\mc{D}}{\mc{Q}}$ is a triangulated functor whose kernel contains $\mc{C}$, then it factors uniquely as $F=Q\circ F_{univ}$ where $\fc{Q}{\mc{D}/\mc{C}}{\mc{Q}}$, see Theorem 2.1.8 \cite{N}. The quotient category $\mc{D}/\mc{C}$ is called Verdier quotient of $\mc{D}$ by $\mc{C}$ and the map $\fc{F_{univ}}{\mc{D}}{\mc{D}/\mc{C}}$ is called Verdier localization map.\\

We denote by $K_0(\mc{D})$ the grothendieck group of a triangulated category $\mc{D}$. By definition, it is the free abelian group generated by isoclasses of indecomposable objects (denoted by $[M]$ for $M$ an object of $\mc{D}$) and they satisfy $[M]+[N]=[P]$ for any three objects that fit into a distinguished triangle $N\rightarrow P\rightarrow M\rightarrow_{+1}$. For an abelian category $\mc{A}$, $K_0(\mc{A})$ is defined in the same way just replacing triangle by short exact sequence. For an additive category $\mc{C}$, we define the split Grothendieck group $K_0(\mc{C}, \oplus)$ with relations given by split short exact sequences. If the category has a monoidal structure, the Grothendieck group becomes a ring under the product $[A][B]=[A\otimes B]$ for any two objects $A,B$ in the cateogry.\\

It is well known that for an abelian (resp. additive) category $\mc{A}$ we have isomorphisms $K_0(\mc{A}) \cong K_0(D^b(\mc{A}))$ (resp. $K_0(\mc{A},\oplus) \cong K_0(K^b(\mc{A}))$) where $D^b(\mc{A})$ (resp. $K^b(\mc{A})$) denotes the bounded derived (resp. homotopy) category of $\mc{A}$. Moreover, If $\mc{D}$ is a triangulated category and $\mc{C}$ is a thick triangulated subcategory, we can compute the Grothendieck group of the Verdier quotient $\mc{D}/\mc{C}$ using the following exact sequence

$$  \xymatrix{ K_0(\mc{C}) \ar[r] & K_0(\mc{D}) \ar[r] & K_0(\mc{D}/\mc{C}) \ar[r] & 0  }    $$

We emphasize that if $\mc{C}$ is not a thick triangulated category the above sequence is not exact in general.\\

We fix the following notation for the rest of the paper: For an abelian monoidal category $\mc{A}$ and a full subcategory $\mc{I}$ of $\mc{A}$ we will denote by $Trg(\mc{I})$ the full triangulated subcategory of $D^b(\mc{A})$ generated by the objects of $\mc{I}$ and we will denote by $\langle \mc{I} \rangle$ the smallest full triangulated subcategory of $D^b(\mc{A})$ that contains $\mc{I}$ and is closed under retracts and tensor products with arbitrary modules.

\subsection{The Beilinson-Bezrukavnikov-Mircovic equivalence $D^b(\mc{U}) \cong K^b(\mc{T})$}

In \cite{BBM} an equivalence between the bounded homotopy category of tilting modules and the bounded derived category of all modules was constructed in the geometric context where modules are replaced by perverse sheaves on a flag manifold. It is easy and well-known how to translate the results of \cite{BBM} into the context of quantum groups. Since we couldn't find a reference in the literature we sketch a proof.

\begin{thm}\label{equivalenceingeneral} The functor $\fc{\gamma}{K^b(\mc{T})}{D^b(\mc{U})}$ induced by inclusion $\T \to \mc{U}$ is an equivalence of triangulated monoidal categories. \end{thm} 

\dem It is obvious that $\gamma$ is a monoidal functor. Since there are no higher extensions between tilting modules and since $K^b(\T)$ is generated by tilting modules as a triangulated category it follows that $\gamma$ is fully faithful.\\

We show that $\gamma$ is essentially surjective. Let $D$ be the subcategory of $D^b(\mc{U})$ (classically) generated by the tilting modules. Since $\gamma$ is fully faithful it suffices to show that $D = D^b(\mc{U})$. For $\lambda\in P^{+}$, let $n(\lambda)$ denote the number of $\mu \in Q^{+}$ such that $\lambda - \mu \in P^{+}$. By \cite{H} Lemma B 13.2 $n(\lambda)< \infty$. We prove by induction on the number $n(\lambda)$ that all the simple modules belongs to $D$. If $n(\lambda) = 0$, then $L(\lambda) = \Delta(\lambda) = T(\lambda) \in D$. Assume $n(\lambda)\neq 0$ and that the result is true for any $\mu\in P^{+}$ such that $n(\mu)<n(\lambda)$. By construction of $T(\lambda)$,  $\Delta(\lambda)$ is a submodule of it and $T(\lambda)$ admits a filtration with sub-quotients $\Delta(\mu)$ for $\mu<\lambda$. Also, $L(\lambda)$ is a quotient of $\Delta(\lambda)$ and $\operatorname{Ker}(\Delta(\lambda) \to L(\lambda))$ admits a filtration with sub-quotients $L(\mu)$, for $\mu<\lambda$. In both cases, $n(\mu) < n(\lambda)$, it follows by induction on the weights that appear in the standard filtration of $T(\lambda)$ that $\Delta(\lambda)\in D$ and so, using the composition series for $\Delta(\lambda)$ we conclude by induction that $L(\lambda) \in D$. Hence $D = D^b(\mc{U})$.
\findem

\begin{cor}\label{homtozero} Any bounded acyclic complex  of tilting modules is contractible. \end{cor}
\findem

The following result was proved in \cite{BKN} using the Balmer spectrum of a triangulated category. Although it is not strictly needed in this paper we opted to include it since it illustrates very well the usefulness of Theorem \ref{equivalenceingeneral}. Given a monoidal category $K$ and an object $M\in K$, let $Tensor_K(M)$ be the thick tensor ideal whose objects are direct summands in $M\otimes N$, $N\in K$.\\

\begin{cor}[\cite{BKN}, 8.2.1 (a)]
Let $\lambda\in P^{+}$. Then $Tensor_{\mc{T}}(T(\lambda)) = Tensor_{\mc{U}}(T(\lambda))\cap \mc{T}$
\end{cor}

\dem Let $V\in \mc{U}$ and assume that $T(\lambda)\otimes V \in Tensor_{\mc{U}}(T(\lambda))\cap \mc{T}$. We must prove that $T(\lambda)\otimes V$ belongs to $Tensor_{\mc{T}}(T(\lambda))$. By Theorem \ref{equivalenceingeneral} there exist a complex $T_V\in K^b(\mc{T})$ such that $T_V\cong V$ in $D^b(\mc{U})$. The complex $T(\lambda)\otimes T_V$ is isomorphic to $T(\lambda)\otimes V$ in $D^b(\mc{U})$. By hypothesis, $T(\lambda)\otimes V\in \mc{T}$; hence, again by Theorem \ref{equivalenceingeneral},  it follows that $T(\lambda)\otimes T_V$ is homotopy equivalent to $T(\lambda)\otimes V$. This fact has the following consequences: First, since $T(\lambda)\otimes T_V$ is bounded and
homotopy equivalent to a complex concentrated in degree $0$ a simple induction shows that $T(\lambda)\otimes T_V$ is homotopy equivalent to a subcomplex which is termwise a direct summand in it and has trivial differential. Second, we conclude that $T(\lambda)\otimes V$ is homotopy equivalent to that subcomplex and since the latter has trivial differential this now implies that the module $T(\lambda)\otimes V$ is isomorphic to a direct summand in its degree $0$ component. Thus $T(\lambda)\otimes V$ is a direct summand in $T(\lambda)\otimes T_V^{0}$. \findem

\subsection{Left triangulated categories and stabilization}\label{lefttriancats}

We recall the notions of a left triangulated category and its stabilization, see \cite{BM}.  A \emph{left triangulated category} $\mc{C}$ consists of an additive category $\mc{C}$, an endofunctor $\Omega: \mc{C} \to \mc{C}$, called the shift functor and a collection of sequences of morphisms
$$\xymatrix{\Omega C \ar[r]^f & A\ar[r]^g & B \ar[r]^h & C }$$
called left triangles (and abbreviated $(A,B,C,f,g,h)$). These data are subject to the axioms LTR1-LTR5 which can be found in \cite{BM}, definition 2.2. Loosely speaking, a left triangulated category satisfies the axioms of a triangulated category with the exception that the shift functor need not be invertible. There is also the notion of triangulated functor between left triangulated categories, this is a functor which sends left triangles to left triangles. Because any triangulated category is left triangulated we can consider triangulated functors between left triangulated and triangulated categories.\\

To any left triangulated category $\mc{C}$ one can associate a triangulated category ${\bf S}(\mc{C})$ called its \emph{stabilization} by formally inverting the shift functor. The objects of ${\bf S}(\mc{C})$ are pairs $(A,n)$ where $A\in \mc{C}$ and $n\in \mb{Z}$, morphisms are given by $\Hom_{{\bf S}(\mc{C})}((A,n),(B,m)) = {\underset{k\geq \max\{n,m\}}{\colim}} \Hom_{\mc{C}}(\Omega^{k-n}(A),\Omega^{k-m}(B))$. The category ${\bf S}(\mc{C})$ is endowed with a triangulated functor $\fc{{\bf S}}{\mc{C}}{{\bf S}(\mc{C})}$ such that for any triangulated category $\mc{D}$ and any triangulated functor $\fc{F}{\mc{C}}{\mc{D}}$, there exist a unique triangulated functor $\fc{F^{\star}}{{\bf S}(\mc{C})}{\mc{D}}$ such that $F^{\star}{\bf S} = F$.

\subsection{Contravariantly finite subcategories}\label{contrfinitesubcat} Let $\mc{A}$ be an abelian category and let $\mc{Y}$ be a full additive subcategory of $\mc{A}$ closed under retracts. A morphism $A \to B$ in $\mc{A}$ is called $\mc{Y}$-epic if any morphism $Y\rightarrow B$, for $Y\in \mc{Y}$, factors through $Y\rightarrow A$. A morphism $\fc{\chi_A}{Y}{A}$ of $\mc{A}$, with $Y\in \mc{Y}$, is called an $\mc{Y}$-\emph{cover} if $\chi_A$ is a $\mc{Y}$-\emph{epic}. The subcategory $\mc{Y}$ is said to be \emph{contravariantly finite} subcategory of $\mc{A}$ if every object of $\mc{A}$ has a $\mc{Y}$-cover. The notion dual to that of $\mc{Y}$-cover is a $\mc{Y}$-\emph{hull} and the notion dual to that of a contravariantly finite subcategory is \emph{covariantly finite} subcategory. If $\mc{Y}$ is both contravariantly and covariantly finite it is called \emph{functorially finite.}\\

An $\mc{Y}$-{\it resolution} of $A\in \mc{A}$ is a complex
\begin{equation}\label{e1}
\cdots \rightarrow Y^1 \rightarrow Y^0 \rightarrow A \rightarrow 0
\end{equation}
where $Y^i\in\mc{Y}$ for $i\geq 0$ and such that the complex
\begin{equation}\label{e2}
\cdots\rightarrow\Hom_{\mc{A}}(Y,Y^1)\rightarrow\Hom_{\mc{A}}(Y,Y^0)\rightarrow\Hom_{\mc{A}}(Y,A)\rightarrow0
\end{equation}
is acyclic for all $Y \in \mc{Y}$. Note that if $\mc{A}$ has enough projectives and the projectives of $\mc{A}$ belongs to $\mc{Y}$ then the exactness of \eqref{e2}
implies that the complex \eqref{e1} is exact as well.\\

We define the {\it contravariant $\mc{Y}$-dimension} of $A$, $\mc{Y}-\dim(A) \in \mb{N} \cup \{\infty\}$ to be the minimal length of a  $\mc{Y}$-resolution of $A$ and
the {\it global contravariant $\mc{Y}$-dimension} of $\mc{A}$,  $$\mc{Y}-gl.\dim(\mc{A}) = \sup\{\mc{Y}-\dim(A) \mid A\in \mc{A}  \}.$$

A complex $(Y,d)$ is said to be $\mc{Y}$-{\it acyclic} in degree $n$ if in the category $\mc{A}$ the morphism $d^{n-1}$ factors as $\xymatrix{ Y^{n-1} \ar[r]^{\epsilon^{n-1}}  & \ker(d^n)\ar[r]  & Y^n  }$ where $\epsilon^{n-1}$ is $\mc{Y}$-epic (here $\ker d^n$ is calculated in $\mc{A}$). $Y$ is called $\mc{Y}$-{\it acyclic} (resp. {\it essentially $\mc{Y}$-acyclic}), if it is $\mc{Y}$-acyclic in degree $n$ for all (resp. for almost all) $n\in \mb{Z}$. We denote by $K^{-,b }(\mc{Y})$ the full subcategory of $K^{-}(\mc{Y})$ consisting of essentially $\mc{Y}$-acyclic complexes.\\

Note that if $\mc{A}$ has enough projectives and these belongs to $\mc{Y}$, then $\mc{Y}$-acyclic in degree $n$ is the same as acyclic in degree $n$ in the usual sense. Then, if $\mc{Y}$ corresponds to the class of projective objects, we say  {\it acyclic} ({\it essentially acyclic}) instead of $\mc{Y}$-acyclic (essentially $\mc{Y}$-acyclic). \\

Assume that $\mc{Y}$ is a contravariantly finite subcategory of $\mc{A}$. For $A,B\in \mc{A}$, denote by $\mc{Y}(A,B)$ the subgroup of $\Hom_{\mc{A}}(A,B)$ consisting of morphisms which factor through an object of $\mc{Y}$. We denote by $\mc{A}/\mc{Y}$  the category whose objects are the same as the objects of $\mc{A}$ and the morphisms are given by $ \Hom_{\mc{A}/\mc{Y}}(A,B) = \Hom_{\mc{A}}(A,B)/\mc{Y}(A,B)$. the category $\mc{A}/\mc{Y}$ is called a \emph{stable category}. By Theorem 2.12 in \cite{BM} and by Theorem 2.2 in \cite{B3}, there exists a left triangulated structure on the stable category $\mc{A}/\mc{Y}$. Moreover, there exists an equivalence of triangulated categories ${\bf S}(\mc{A}/\mc{Y}) \cong K^{-,b}(\mc{Y})/K^b(\mc{Y})$. In particular, if the category $\mc{A}$ has enough projectives and $\mc{P}$ denotes the full subcategory of projectives in $\mc{A}$, then ${\bf S}(\mc{A}/\mc{P}) \cong D^b(\mc{A})/K^b(\mc{P})$. Furthermore, if $\mc{A}$ is a Frobenius category, ${\bf S}(\mc{A}/\mc{P})\cong\mc{A}/\mc{P}$, see \cite{B} Theorem 3.8 and Corollary 3.9 for details.

\section{Triangulated categories associated  to the big quantum group}

In this section we study the Verdier quotient $K^b(\mc{T})/K^b(\mc{N})$ as a triangulated analogue of the category $\mc{F}$. We show that it is generated as a triangulated category by the simple modules with highest weights in the principal alcove and that its Grothendieck ring coincides with the fusion ring $\mc{R}$. We also investigate the category $\mc{N}_{\mc{U}}$ consisting of all negligible modules in the category $\mc{U}$, the Verdier quotient $D^b(\mc{U})/\langle\mc{N}_{\mc{U}}\rangle$ and its Grothendieck ring.

\subsection{The category $K^b(\mc{T})/K^b(\mc{N})$ and its Grothendieck ring}

Let $\fc{\pi}{K^b(\mc{T})}{K^b(\mc{F})}$ be the functor induced by the functor defined in equation \eqref{TtoFpart}. Then, by definition $\pi(X) = X^{\mc{F}}$. Note that $\pi$ is a monoidal functor because $\mc{N}$ is a tensor ideal in $\mc{T}$. Let $\fc{\epsilon}{K^b(\mc{F})}{K^b(\mc{T})}$ be the natural functor induced by the inclusion $\mc{F}\hookrightarrow \mc{T}$.\\

\begin{prop}\label{closedretracts} The category $K^b(\mc{N})$ is closed under retracts. \end{prop}

\dem Let $Y \in K^b(\mc{T})$ be a retract of $X \in K^b(\mc{N})$,
that is, there are maps $i: Y \to X$ and $r:X \to Y$ such that $ri$ is homotopic to $Id_Y$. We must show that $Y$ is homotopy equivalent to an object of $K^b(\mc{N})$.

For each $n\in\mb{Z}$ we have that $Y^{n} = Y^{n}_F\oplus Y^{n}_N$,
where $Y^{n}_F \cong (Y^n)^{\mc{F}}$ and $Y^{n}_N\in \mc{N}$. If
$Y^{n}_F = 0$ for all $n$ we are done. If not, let $n$ be the
smallest index such that $Y^{n}_F \neq 0$.

Let $d$ denote the differential in $Y$ and write $d^n = \begin{pmatrix} d_{11}^n & d_{12}^n \\ d_{21}^n &
d_{22}^n \end{pmatrix}$. Note that $d_{11}^{i+1}d_{11}^{i}=0$ for
$i\geq n$, because $d^{i+1}d^{i}=0$ and its $(1,1)$-component has
the form $d_{11}^{i+1}d_{11}^{i} +  d_{12}^{i+1}d_{21}^{i} =0$,
but $d_{12}^{i+1}d_{21}^{i}=0$ since only the zero map between
fusion modules can factor through a negligible tilting module.
Using the homotopy between $ri$ and $Id_{Y}$ it is easy to see
that $\fc{d_{11}^n}{Y_F^n}{Y_F^{n+1}}$ is a split monomorphism.
Choose a submodule $(Y_F^n)^{\perp} \subseteq Y_F^{n+1}$ such that
$Y_F^{n+1} = \operatorname{Im} d_{11}^n \oplus
(Y_F^n)^{\perp}$. %and the morphism $d_{11}^n$ is the identity restricted to $Y_F^{n-1}$.
For an element $x\in Y_F^{n+1}$ we write $x=x' + x^{\perp}$ for
its components with respect to this decomposition. Define a
complex $\tilde{Y} \in K^b(\mc{T})$ as follows
$$\xymatrix{ \tilde{Y}: \cdots \ar[r] & Y_N^{n-1} \ar[r]^{\partial^{n-1}} & Y_N^n \ar[r]^(0.3){\partial^{n}}
& (Y_F^n)^{\perp}\oplus Y_N^{n+1} \ar[r]^{\partial^{n+1}} &
Y_F^{n+2}\oplus  Y_N^{n+2} \ar[r] & \cdots } $$ where the
differential $\partial$ is defined as follows:  $\partial^{n-1} =
d_{22}^{n-1}$, $\partial^{n}(y)=((d_{12}^ny)^{\perp}, d_{22}^ny -
d_{21}^n((d_{12}^nx)_n))$ and $\partial^j = d^j$, for $j \neq
n,n-1$.

In order to verify that $\partial\circ \partial=0$, consider the sequence of maps $f^i:Y^i\to\tilde{Y}^i$ defined by $f^i(x,y)=(x,y)$ for $i\neq n , n+1$, $f^n(x,y)=y$ and $f^{n+1}(x,y)=(x^{\perp}, y-d_{21}^{n}x')$. Then, $f^{i+1} d^i = \partial^{i}f_i$ for all $i$ and since each $f^i$ is surjective we conclude that $\partial \circ \partial = 0$. Hence, $f = \{f^i\}$ is a morphism of complexes.

We now show that $\tilde{Y}$ is homotopy equivalent to $Y$. Define a
morphism $g:\tilde{Y}\to Y$ as follows: $g^{i}=Id_{Y^i}$ for
$i\leq n-1$, $g^{n}(x)=(-(d_{12}^nx)',x)$ and $g^{i}$ are the
inclusions for $i\geq n+1$. It follows that $gf$ is homotopic to $Id_{\tilde{Y}}$ and $fg$ is homotopic to $Id_{Y}$.

Iterating this process we will after a finite number of steps
obtain a complex in $K^b(\mc{N})$ which is homotopic to $Y$ and we have thus proved that $K^b(\mc{N})$ is closed under retracts.

\findem

\begin{prop}\label{newpartprop} The functor $\pi$ factors through a monoidal functor $\fc{\overline{\pi}}{K^b(\mc{T})/K^b(\mc{N})}{K^b(\mc{F})}$. Moreover, $\pi \circ \epsilon \cong Id_{K^b(\mc{F})}$.  \end{prop}

\dem The first part holds since $\pi(K^b(\mc{N})) = 0$. The final part is evident. \findem

\begin{prop}\label{noadjointpi} The functor $\pi$ does not admit neither a right nor a left adjoint. \end{prop}

\dem  Suppose $\fc{r}{K^b(\mc{F})}{K^b(\mc{T})}$ is a right adjoint of the functor $\pi$. Then, $\Hom_{K^b(\mc{T})}(N, rF)$ $\cong \Hom_{K^b(\mc{F})}(\pi N,F)=0$ for all $N\in K^b(\mc{N})$ and $F\in K^b(\mc{F})$, because $\pi N=N^{\mc{F}}\cong 0$. Since the projective objects of $\mc{U}$ belongs to $\mc{N}$ and $\mc{U}$ has enough projectives, we get that $rF$ is acyclic. But any exact complex in $K^b(\mc{T})$ is necessarily contractible by Corollary \ref{homtozero}. Hence $rF\cong 0$. Now, if we take $F=\mb{C}[0]$,
$$\mb{C} \cong \Hom_{K^b(\mc{F})}(\pi\mb{C}[0],\mb{C}[0]) \cong \Hom_{K^b(\mc{T})}(\mb{C}[0],r\mb{C}[0]) \cong \Hom_{K^b(\mc{T})}(\mb{C}[0],0)\cong 0,$$ which is a contradiction. Hence, the functor $\pi$ does not admits a right adjoint. Similarly, there is no left adjoint because all injective modules belong to $\mc{N}$. \findem

\begin{prop}\label{generatorquotient} The category $K^b(\mc{T})/K^b(\mc{N})$ is generated as a triangulated category by the classes of the objects $L(\lambda)$ for $\lambda\in C_{\ell}$. \end{prop}

\dem By definition $K^b(\mc{T})$ and therefore its quotient $K^b(\mc{T})/K^b(\mc{N})$ are generated by $T(\lambda)$, for $\lambda \in P^{+}$. For $\lambda \notin C_{\ell}$ we have $T(\lambda) = 0$ in $K^b(\mc{T})/K^b(\mc{N})$ and for $\lambda \in C_{\ell}$ we  have $T(\lambda) = L(\lambda)$. \findem

\begin{prop}\label{TmodNnotss} The category $K^b(\mc{T})/K^b(\mc{N})$ is not semisimple. \end{prop}

\dem If the category $K^b(\mc{T})/K^b(\mc{N})$ were to be semisimple, then by the previous proposition the functor $\fc{\overline{\pi}}{K^b(\mc{T})/K^b(\mc{N})}{K^b(\mc{F})}$ of Proposition \ref{newpartprop} would be an equivalence and so the functor $\fc{\pi}{K^b(\mc{T})}{K^b(\mc{F})}$ will admits an adjoint, which is not possible by Proposition \ref{noadjointpi}. \findem

We proceed now to compute the Grothendieck ring of the quotient category $K^b(\mc{T})/K^b(\mc{N})$. By \cite{AP} 3.19, we have that the Grothendieck ring of the category $\mc{U}$ has $[T(\lambda)]$, $\lambda\in P^{+}$ as a $\mb{Z}$-basis. It has decomposition of the form $K_0(\mc{U}) \cong \mc{R}\oplus \mc{R}^{+}$ where $ \mc{R} \cong span_{\mb{Z}}\{ [T(\lambda)] | \lambda\in C_{\ell} \} $ and $ \mc{R}^{+} \cong span_{\mb{Z}}\{ [T(\lambda)] | \lambda \in P^{+}\setminus C_{\ell} \} $. Moreover, $K_0(\mc{N}, \oplus)\cong\mc{R}^{+}$ and as we defined earlier $K_0(\mc{F}) = \mc{R}$.\\

Define a ring homomorphism $\fc{\varphi}{K_0(K^b(\mc{T})/K^b(\mc{N}))}{\mc{R}}$ by $[L(\lambda) \mod K^b(\mc{N})]\mapsto[L(\lambda)]$, where $\lambda\in C_{\ell}$.

\begin{prop}\label{Grothringquotient01} $\varphi$ is an isomorphism. \end{prop}

\dem At the level of Grothendieck groups we have exact sequence $K_0(K^b(\N^\lambda))\rightarrow K_0(K^b(\T^\lambda)) \rightarrow K_0(K^b(\T^\lambda)/K^b(\N^\lambda))\rightarrow 0$. Because $K^b(\N^\lambda)$ is closed under retracts the Euler characteristic of the quantum trace vanishes on $K^b(\N^\lambda)$, and we have well defined map $K_0(K^b(\T^\lambda)/K^b(\N^\lambda))\rightarrow \mb{C}$. By Proposition \ref{generatorquotient} the category $K^b(\mc{T})/K^b(\mc{N})$ is generated by the simple modules with weights in the principal alcove. In particula, for $\lambda \in C_{\ell}$, $K^b(\mc{T}^\lambda)/K^b(\mc{N}^\lambda)$ is generated by the simple module of weight $\lambda$, which has non-zero quantum dimension, and hence $K^b(\mc{T}^\lambda)/K^b(\mc{N}^\lambda)$ is a free abelian group of rank one. Then $\varphi$ defines map from $K_0(K^b(\T)/K^b(\N)$ to $K^b(\F)$ which send generators to generators. Because both are free abelian groups of the same rank, $\varphi$ is an isomorphism. \findem

\subsection{The category $D^b(\mc{U})$ and a characterization of $\langle \mc{N} \rangle$}

Recall the block decomposition $\mc{U}=\oplus_{\lambda\in P/W_{\ell}}\mc{U}^{\lambda}$ given by the linkage principle. We say that a block $\mc{U}^{\lambda}$ is singular if $\lambda$ is a singular weight. Denote by $\langle D^b(\mc{U})_{sing} \rangle$ the smallest triangulated subcategory of $D^b(\mc{U})$ that contains the singular blocks of $\mc{U}$ and is closed under retracts and tensor products with arbitrary modules. Equivalently,  $\langle D^b(\mc{U})_{sing} \rangle$ is the smallest triangulated subcategory closed under retracts and tensor ideal of $D^b(\mc{U})$ which contain $L(\lambda)$ for $\lambda \in P^{+}$ singular.\\

Recall that $Trg(\mc{N})$ denotes the triangulated subcategory of $D^b(\mc{U})$ generated by $\mc{N}$. Thus $Trg(\mc{N}) = \gamma(K^b(\mc{N}))$ and hence $Trg(\mc{N})$ is a tensor ideal. By Proposition \ref{closedretracts}, $K^b(\mc{N})$ is closed under retracts, so $\langle\mc{N}\rangle = Trg(\mc{N})$ (recall that $\langle\mc{N}\rangle$ is the smallest full triangulated subcategory of $D^b(\mc{U})$ that contains $\mc{N}$ and is closed under retracts and tensor products with arbitrary modules, see section \ref{GTI}).

\begin{thm}\label{NidealDU} $\langle D^b(\mc{U})_{sing}\rangle = \langle \mc{N} \rangle$. \end{thm}

\dem As a triangulated category and tensor ideal $\langle D^b(\mc{U})_{sing} \rangle$ is generated by $M\in \mc{U}^{\lambda}$, $\lambda$ singular.  By  Theorem \ref{equivalenceingeneral} and the linkage principle $M\cong T$ in $D^b(\mc{U})$ for some
$T \in K^b(\mc{T}^{\lambda})$. Since, $\mc{T}^\lambda \subset \mc{N}$ we have proved that $\langle D^b(\mc{U})_{sing} \rangle \subseteq \langle \mc{N} \rangle$.

On the other hand, for any $T(\mu)\in \mc{N}$, we have that $\mu\in \overline{A}$, where $A$ is an alcove different form $C_{\ell}$. If $\mu\in \overline{A}\setminus A$, $\mu$ is singular and we are done. If $\mu\in A$ is regular, then there exists $\lambda \in \overline{A}\setminus A$ such that $\mu-\lambda\in P^{+}$. Put $E=T(\mu-\lambda)$. Then, $T(\mu)$ is a direct summand of $T(\lambda)\otimes E$. Since $\lambda$ is singular,  $T(\lambda)\in \mc{U}^{\lambda}$ and so $T(\lambda)\otimes E$ belongs to $\langle D^b(\mc{U})_{sing} \rangle$. Therefore, $T(\mu)\in \langle D^b(\mc{U})_{sing} \rangle$ because the ideal $\langle D^b(\mc{U})_{sing} \rangle$ is closed under retracts. Thus $\langle \mc{N} \rangle$ is contained in $\langle D^b(\mc{U})_{sing} \rangle$. \findem

Theorem \ref{equivalenceingeneral} implies that the category $K^b(\mc{T})/K^b(\mc{N})$ is equivalent to the category $D^b(\mc{U})/\langle \mc{N}\rangle$. Thus, by Proposition \ref{Grothringquotient01} we have the following

\begin{cor}\label{Grothringquotient} The ring $K_0(D^b(\mc{U})/\langle \mc{N} \rangle)$ is isomorphic to the ring $\mc{R}$. \end{cor}
\findem

\subsection{The category $D^b(\mc{U})/\langle\mc{N}_{\mc{U}}\rangle$ and its Grothendieck ring}\label{conj}

Recall that $\mc{N}_{\mc{U}}$ is the full subcategory of $\mc{U}$ consisting of all negligible modules, that is, modules such that the quantum trace of each of its endomorphisms vanishes.

\begin{lem}[see \cite{EOSS}]\label{lemnegind} \begin{itemize}
\item[1)] $\N_{\mc{U}}$ is closed under direct sums and summands.
\item[2)] Let $M\in \mc{U}$ be indecomposable, then $M\in \N_{\mc{U}}$ if and only if $\dim_q(M)=0$.  
\end{itemize} \end{lem}

\dem $1)$ Follows because the decomposition of any module into indecomposable modules has finitely many components.\\
 $2)$ On one side, by definition the trace of the identity morphism of an indecomposable negligible module is zero. On the other side, if $f$ is an endomorphism of $M$ it can be written as a sum of a scalar and a nilpotent morphisms $f_e$ and $f_n$ respectively, because $M$ is indecomposable. The trace of the former is zero because $\dim_q(M)=\Tr_q(1_M)=0$ and the trace of the latter is zero since $f_nK_{2\rho}$ a nilpotent operator and $M$ a direct sum of its weight spaces. \findem

Consider the full triangulated subcategory and tensor ideal $\langle \mc{N}_{\mc{U}} \rangle$ of $D^b(\mc{U})$. We want to study the Verdier quotient $D^b(\mc{U})/\langle \mc{N}_{\mc{U}} \rangle$.\\

Let $\lambda$ be a dominant weight, denote by $\mc{N}^{\lambda}=\mc{N}\cap\mc{U}^{\lambda}$ and $\mc{N}_{\mc{U}}^{\lambda}=\mc{N}_{\mc{U}}\cap\mc{U}^{\lambda}$ the respective blocks. Since indecomposable modules has composition factors with linked highest weights and there are no morphisms between simples unless the highest weights are the same, by Lemma \ref{lemnegind} we have

$$ \mc{N}_{\mc{U}} \cong \bigoplus_{\lambda\in P/W_{\ell}} \mc{N}_{\mc{U}}^{\lambda}, \quad \quad K^b(\mc{N}_{\mc{U}}) \cong \bigoplus_{\lambda\in P/W_{\ell}} K^b(\mc{N}_{\mc{U}}^{\lambda})$$

and from this we conclude that $$ D^b(\mc{U})/\langle\mc{N}_{\mc{U}}\rangle \cong \bigoplus_{\lambda\in P/W_{\ell}} D^b(\mc{U}^{\lambda})/\langle\mc{N}_{\mc{U}}^{\lambda}\rangle.$$

Consider the Euler characteristic of the quantum trace, that is, the morphisms $\Tr_q^{\bullet} : D^b(\mc{U}) \rightarrow \mb{C}$ defined by $\Tr_q^{\bullet}(X)=\sum_i (-1)^{i}\Tr_q(X^{i})$, where $X^{i}$ are the components of the complex $X$. This is a well defined map because the sum is finite. Let $[\Tr_q] : K_0(D^b(\mc{U})) \rightarrow \mb{C}$ the induce trace map in the Grothendieck group.\\

%Since the tensor ideals $\langle\mc{N}_{\mc{U}}^{\lambda}\rangle$ are closed under retracts, 

\noindent If $[\Tr_q]([N])=0$ for any object $N$ in $\langle\mc{N}_{\mc{U}}^{\lambda}\rangle$, using the exact sequence (see section \ref{GTI})

$$  \xymatrix{ K_0(\langle\mc{N}_{\mc{U}}^{\lambda}\rangle) \ar[r] & K_0(D^b(\mc{U}^{\lambda})) \ar[r] & K_0(D^b(\mc{U}^{\lambda})/\langle\mc{N}_{\mc{U}}^{\lambda}\rangle) \ar[r] & 0  }    $$

\noindent we obtain a well defined map $[\Tr_q]: K_0(D^b(\mc{U}^{\lambda})/\langle\mc{N}_{\mc{U}}^{\lambda}\rangle) \rightarrow \mb{C}$. Then, using that $ \mc{N}_{\mc{U}} \cong \bigoplus_{\lambda\in P/W_{\ell}} \mc{N}_{\mc{U}}^{\lambda}$ and the short exact sequence 

$$  \xymatrix{ K_0(\langle\mc{N}_{\mc{U}}\rangle) \ar[r] & K_0(D^b(\mc{U})) \ar[r] & K_0(D^b(\mc{U})/\langle\mc{N}_{\mc{U}}\rangle) \ar[r] & 0  }    $$

we have well defined map $[\Tr_q]: K_0(D^b(\mc{U})/\langle\mc{N}_{\mc{U}}\rangle) \rightarrow \mb{C}$.\\

We still do not known if this map exists or not. Moreover, we do not know if $\langle\N_{\mc{U}}\rangle$ is the whole derived category $D^b(\mc{U})$ or weather $\langle\N_{\mc{U}}\rangle = \langle\N\rangle$. One way to know if the map $[\Tr_q]: K_0(D^b(\mc{U})/\langle\mc{N}_{\mc{U}}\rangle) \rightarrow \mb{C}$ exists is to prove that $Trg(\N_{\mc{U}})\cong\langle\N_{\mc{U}}\rangle$, which means that all the retracts in $\langle\N_{\mc{U}}\rangle$ are negligible, so $[\Tr_q](\langle\N_{\mc{U}}\rangle)=0$ and the map $[\Tr_q]$ is defined. Moreover, $\langle\N_{\mc{U}}\rangle$ becomes a proper ideal. Using this and that each of the categories $D^b(\mc{U}^{\lambda})/\langle\mc{N}_{\mc{U}}^{\lambda}\rangle$ are generated by the irreducible module with highest weight $\lambda$, we are able to prove that $K_0(D^b(\mc{U}^{\lambda})/\langle\mc{N}_{\mc{U}}^{\lambda}\rangle)\cong \mb{Z}$ and so $K_0(D^b(\mc{U})/\langle\mc{N}_{\mc{U}}\rangle)\cong \mc{R}$.\\  

We will investigate in Proposition \ref{closedretractsequiv} when $Trg(\N_{\mc{U}})\cong\langle\N_{\mc{U}}\rangle$ holds.

\section{The stable category ${\bf S}(\mc{U}/\mc{N})$}

In this section, we show that the category $\mc{N}$ is functorially
finite in the category $\mc{U}$. This allow us to define the
stable category ${\bf S}(\mc{U}/\mc{N})$ and construct a quotient
functor from this category to the category
$D^b(\mc{U})/\langle\mc{N}\rangle$. As a consequence, we conclude
that the Grothendieck ring of the category ${\bf
S}(\mc{U}/\mc{N})$ is an enhancement for the fusion ring $\mc{R}$.

\subsection{$\mc{N}$ is functorially finite in $\mc{U}$}

We start by showing that the category $\mc{N}$ is functorially
finite in $\mc{U}$, (see Section \ref{contrfinitesubcat} for definitions).\\

Given a tilting module $T$, a Weyl module $\Delta$ and a simple module $L$, denote by $[T:\Delta]$ the standard flag multiplicity of $\Delta$ in $T$ and by $[\Delta:L]$ the Jordan - H{\"o}lder multiplicity of $L$ in $\Delta$.

\begin{lem}\label{finsimpletilting} Let $\mu\in P^{+}$ and $L(\mu)$ be the simple $\U_q$-module of weight $\mu$. Then,
there exists only finitely many $\lambda\in P^{+}$ such that $\Hom_{\mc{U}}(T(\lambda),L(\mu))\neq 0$.
Hence, for any $V\in \mc{U}$, there exists finitely many $\lambda\in P^{+}$ such that $\Hom_{\mc{U}}(T(\lambda),V)\neq 0$.
\end{lem}

\dem Fix $L(\mu)$ for $\mu\in P^{+}$. Consider $\lambda\in P^{+}$ and the indecomposable tilting module $T(\lambda)$. Then

 $$\begin{array}{c l}
\dim_{\mb{C}}(\Hom_{\mc{U}}(T(\lambda),L(\mu))) & = [T(\lambda): L(\mu)] \\
& = \sum_{\nu\in P^{+}} [T(\lambda) : \Delta(\nu)][\Delta(\nu): L(\mu)] \\
& = \sum_{\nu\in P^{+}} [T(\lambda) : \Delta(\nu)][T(\overline{\mu}): \Delta(\nu)] \\
& = \sum_{\nu\in P^{+}} n_{\nu\lambda}(1)n_{\nu\overline{\mu}}(1) \\
\end{array} $$
\\
The second equality is given by BGG-reciprocity. $T(\overline{\mu})$ denotes the projective cover of $L(\mu)$ and
$n_{\nu\lambda}$, $n_{\nu\overline{\mu}}$ denotes the parabolic
Kazhdan-Lusztig polynomials, see \cite{S}. If $n_{\nu\lambda}(1)n_{\nu\overline{\mu}}(1)\neq 0$ then $\nu\leq\overline{\mu}$ and $\nu\leq \lambda$. If $[\Delta(\nu):L(\mu)]\neq 0$ then $\mu\leq\nu$. Therefore, $\mu\leq \nu\leq\overline{\mu}$. Hence, there are only finitely many tilting modules $T(\lambda)$ with standard composition factors having weights $\nu$ such that $\mu\leq \nu\leq\overline{\mu}$.

The last statement of the lemma follows because the category $\mc{U}$ is a finite length category. \findem

\begin{thm}\label{Nresolutions} The category $\mc{N}$ is functorially finite in $\mc{U}$. \end{thm}

\dem By definition, it is sufficient to show that  any $V\in \mc{U}$ admits an $\mc{N}$-cover. For each $\lambda \in P^{+}\setminus C_{\ell}$ let $n_\lambda=\dim\Hom_{\mc{U}}(T(\lambda), V)$. Then each $n_{\lambda}$ is finite and it is zero for almost all $\lambda$ by Lemma \ref{finsimpletilting}. Let $N_V := \bigoplus_{\lambda\in P^{+}\setminus C_{\ell}} T({\lambda})^{n_\lambda}\in \mc{N}$. Let $\fc{can}{N_V}{V}$ be the canonical map. Then by construction any map $N\rightarrow V$, for $N\in \mc{N}$ factors through $can$. Also note that $can$ is surjective since $\mc{U}$ has enough projectives and the projectives belongs to $\mc{N}$. Thus, $can$ is an $\mc{N}$-cover. \findem

Theorem \ref{Nresolutions} implies that any object $A\in \mc{U}$ admits an
$\mc{N}$-resolution $N_A \to A$, where $N_A \in K^{-}(\mc{N})$
lives in non-positive degrees. Since $\mc{N}$ contains all the
projectives and $\mc{U}$ has enough projectives it follows that
the complex $N_A \to A$ is automatically acyclic in the usual
sense. Recall that $\mc{N}-gl.\dim(\mc{U})$ is the supremum of all the minimal lengths of $\mc{N}$-resolutions for objects in $\mc{U}$.

\begin{prop}\label{globalNdim} $\mc{N}-gl.\dim(\mc{U})$ is infinite. \end{prop}

\dem Let $F\in \mc{T}\setminus \mc{P}$, $P\twoheadrightarrow F$ be
a projective cover and let $A=\ker(P\twoheadrightarrow F)$. We
show that $A$ cannot admit a finite $\mc{N}$-resolution. If it did, we would have a finite acyclic complex $N_A \to A$ where $N_A \in
K^{b}(\mc{N})$. But this would give the acyclic complex $N_A \to P
\to F$ in $K^b(\mc{T})$. Hence, by Corollary \ref{homtozero}, the
latter is contractible, which is impossible since $P
\twoheadrightarrow F$ is not split. \findem

\subsection{The category $K^{-,b}(\mc{N})/K^b(\mc{N})$}

Theorem \ref{Nresolutions} allow us to construct the left triangulated category $\mc{U}/\mc{N}$ and its stabilization ${\bf S}(\mc{U}/\mc{N})$, which is a triangulated category equivalent to the Verdier quotient category $K^{-,b}(\mc{N})/K^b(\mc{N})$, see sections \ref{lefttriancats} and \ref{contrfinitesubcat}.\\

Recall that the objects of the category $K^{-,b}(\mc{N})$ are essentially $\mc{N}$-acyclic complexes. Since projective modules are negligible, the objects of $K^{-,b}(\mc{N})$ are essentially acyclic, that is, they have bounded cohomologies. For that reason, the image of the composition $K^{-,b}(\mc{N})\rightarrow K^{-}(\mc{N}) \rightarrow D^{-}(\mc{U})$ is contained in $D^b(\mc{U})$ and defines a triangulated functor $\fc{F}{K^{-,b}(\mc{N})}{D^b(\mc{U})}$, which is the identity on objects and induces a functor
$$ \fc{\overline{F}}{K^{-,b}(\mc{N})/K^b(\mc{N})}{D^b(\mc{U})/\langle \mc{N} \rangle} .$$

Theorem \ref{Nresolutions} asserts that any object of $\mc{U}$ has an $\mc{N}$-cover and by construction this cover is an epimorphism. In particular, tilting modules have $\N$-covers. Since the category of negligible tilting modules is closed under taking direct sums and under isomorphisms, for any complex $X$ in $D^b(\mc{U})$ there exist a complex $N_X$ in $K^{-,b}(\mc{N})$ and a morphism $N_X\to X$ which is a quasi-isomorphism, Theorem I.7.5 of \cite{Betal}. Note that in the case $X$ is a module in $\mc{U}$, then $N_X$ is an $\mc{N}$-resolution of $X$. We define the functor $\fc{G}{K^b(\mc{T})}{K^{-,b}(\mc{N})}$ by $G(X)=N_X$. Then $G$ induces a functor $$\fc{\overline{G}}{K^b(\mc{T})/K^b(\mc{N})}{K^{-,b}(\mc{N})/K^b(\mc{N})}$$

\begin{lem}\label{fullyfaithfuliota} We have the following natural isomorphisms for $F \in \mc{F}$ and $N\in \mc{N}$:
\begin{enumerate}
\item[1)] $\Hom_{K^{-,b}(\mc{N})}(G(F),N[0]) \cong \Hom_{K^b(\mc{T})}(F[0],N[0])$ and $\Hom_{K^{-,b}(\mc{N})}(G(F),N[i])\cong 0$ for $i\neq 0$.
\item[2)] $\Hom_{K^{-,b}(\mc{N})}(N[0], G(F)) \cong \Hom_{K^b(\mc{T})}(N[0],F[0]) $
\item[3)] $\Hom_{K^{-,b}(\mc{N})}(G(F), G(F)[0]) \cong \Hom_{K^b(\mc{T})}(F[0],F[0]) $ and\\
 $\Hom_{K^{-,b}(\mc{N})}(G(F), G(F)[i])\cong 0$ for $i\neq 0$.
\end{enumerate}
\end{lem}

\dem Recall that $\Ext^{i}_{\mc{U}}(T,T')\cong 0$ for $i>0$ and any two tilting modules $T$ and $T'$. $1).$ Because $G(F)$ is a deleted $\mc{N}$-resolution of $F$, it satisfies $H^0(G(F)) \cong F$ and the other cohomologies are zero. So, we have
$$\Hom_{K^{-,b}(\mc{N})}(G(F), N[i]) \cong \Ext^i_{\mc{U}}(H^{0}(G(F)), N) \cong \Ext^i_{\mc{U}}(F, N)$$ 
which is isomorphic to zero if $i\neq 0$ and is isomorphic to $\Hom_{\mc{U}}(F,N)$ if $i=0$. then 
$$\Hom_{K^{-,b}(\mc{N})}(G(F), N[0]) \cong\Hom_{\mc{U}}(F,N) \cong \Hom_{K^b(\mc{T})}(F[0],N[0]).$$ 
$2).$ Since $G(F) = N_F$ is an $\mc{N}$-resolution of $F$ (which components will be denoted by $N^{i}_F$), the following complex is acyclic
$$\xymatrix{\cdots \ar[r] & \Hom_{\mc{N}}(N, N_F^2) \ar[r] & \Hom_{\mc{N}}(N, N_F^1) \ar[r] & \Hom_{\mc{N}}(N, F) \ar[r] & 0}.$$
Therefore, $$ \Hom_{\mc{N}}(N,F) \cong \Hom_{\mc{N}}(N,N_F^1) / \im(\Hom_{\mc{N}}(N,N_F^2)\rightarrow \Hom_{\mc{N}}(N,N_F^1)) $$ but the right hand side is $\Hom_{K^{-,b}(\mc{N})}(N[0], G(F))$ and the left hand side is $\Hom_{K^b(\mc{T})}(N[0],F[0])$. \\
$3).$ $\Hom_{K^{-,b}(\mc{N})}(G(F), G(F)) \cong \Hom_{D^b(\mc{U})}(G(F), H^{0}(G(F))) \cong \Hom_{\mc{U}}(F,F)$\\ $\cong \Hom_{K^b(\mc{T})}(F[0],F[0]) $. \findem

We denote by $K^{-,b,ex}(\mc{N}) \subset K^{-}(\mc{N})$ the subcategory of acyclic complexes. Note that  $K^{-,b,ex}(\mc{N})=\ker(F)$. We define a triangulated functor

$$  \fc{F'}{(K^{-,b}(\mc{N})/K^b(\mc{N}))/K^{-,b,ex}(\mc{N})}{D^b(\mc{U})/\langle \mc{N} \rangle} $$

and a the triangulated functor $$ G': D^b(\mc{U})/\langle \mc{N} \rangle \rightarrow (K^{-,b}(\mc{N})/K^b(\mc{N}))/K^{-,b,ex}(\mc{N}) $$ by $G':= Q'\circ G \circ \overline{\gamma}^{-1}$, where the functor $Q':K^{-,b}(\mc{N})/K^b(\mc{N}) \rightarrow (K^{-,b}(\mc{N})/K^b(\mc{N}))/K^{-,b,ex}(\mc{N})$ and $\overline{\gamma}:K^b(\mc{T})/K^b(\mc{N}) \rightarrow D^b(\mc{U})/\langle\mc{N}\rangle$ is the equivalence induced by $\gamma$

\begin{thm}\label{thetheorem} The functor $F'$ is an equivalence of triangulated categories. \end{thm}

\dem By Lemma $\ref{fullyfaithfuliota}$ the functor $\overline{G}$ is fully faithful. If $X$ belongs to $D^b(\mc{U})/\langle \mc{N} \rangle$, then $G(X) = N_X$ and $N_X\to X$ is an isomorphism. Hence, $\overline{F}\circ\overline{G}(X) \to \overline{F}(X) =X$ is an isomorphism. If $X$ belongs to $K^{-,b}(\mc{N})/K^b(\mc{N})$, there exist triangle $ \overline{G} \circ\overline{F}(X) \to X\to C \to_{+1}$ for some $C$. Applying $\overline{F}$ to this triangle and using the fact that $\overline{F}\circ\overline{G}\circ\overline{F}(X) = \overline{F}\circ\overline{G}(\overline{F}(X)) \cong \overline{F}(X)$ we get triangle $\overline{F}(X) \to \overline{F}(X)\to \overline{F}(C) \to_{+1}$.  Therefore, $\overline{F}(C) \cong  0$, and $C\in \ker(\overline{F})$. Then, in $(K^{-,b}(\mc{N})/K^b(\mc{N}))/K^{-,b,ex}(\mc{N})$, we get that $C\cong 0$. So, $\overline{G}\circ\overline{F}(X) \to X$ is an isomorphism. Hence, $G'$ is an inverse for $F'$, hence $F'$ is an equivalence of categories. \findem

We can relate the above quotient categories with the stabilization of the left homotopy pair $(\mc{U}, \mc{N})$. Denote by $\mathscr{K}$ the image of $K^{-,b,ex}(\mc{N})/K^b(\mc{N})$ under the equivalence of ${\bf S}(\mc{U}/\mc{N})$ with $K^{-,b}(\mc{N})/K^b(\mc{N})$.

\begin{cor} ${\bf S}(\mc{U}/\mc{N})/\mathscr{K} \rightarrow K^b(\mc{T})/K^b(\mc{N}) \rightarrow D^b(\mc{U})/\langle \mc{N} \rangle$ are equivalences of triangulated categories. \end{cor}

\dem By Theorem \ref{Nresolutions} and section \ref{contrfinitesubcat}, ${\bf S}(\mc{U}/\mc{N})$ is equivalent to  $K^{-,b}(\mc{N})/K^b(\mc{N})$. This result plus the Theorems \ref{thetheorem} and \ref{equivalenceingeneral} prove the desired equivalences. \findem

\subsection{Grothendieck rings}

We have not been able to explicitly describe $K_0(K^{-,b}(\mc{N})/K^b(\mc{N}))$. However we shall see that we have a surjective and non-injective, unless $\g=\mf{sl}_2$, map to $\mc{R}$. For this purpose, for $N\in \mc{N}$ define the functions $\fc{a_N,b_N}{K^{-,b}(\mc{N})}{\mb{Z}}$ by

$$ a_N(X) = \sum_{i\in\mb{Z}} (-1)^{i}\dim\Hom_{K^{-}(\mc{N})}(X,N[i]) $$
$$ b_N(X) = \sum_{i\in\mb{Z}} (-1)^{i}\dim\Hom_{K^{-}(\mc{N})}(N[0],X[i]) $$

\begin{lem}\label{welldefiab} For every $N\in\mc{N}$, the functions $a_N$ and $b_N$ are well-defined (i.e. their defining sums are finite) and additive with respect to distinguished triangles.\end{lem}

\dem Let $X$ be a complex in $K^{-,b}(\mc{N})$. Let $r_E$ be the smallest integer such that $X^i$ is $\mc{N}$-acyclic for $i<r_E$. Because the complex $X$ is bounded above, let $r_X$ be the smallest integer such that $X^{i}=0$ for any $i>r_X$.\\

We verify that the function $b_N$ is well defined. If $i<<0$, then $\Hom_{K^{-}(\mc{N})}(N[0],X[i])=0$. On the other hand, if we pick $i>>0$, since $X$ is bounded above, $\Hom_{K^{-}(\mc{N})}(N[0],X[i])=0$.\\

Let see that the function $a_N$ is well-defined. Because $X$ is bounded above, for any $i<-|r_X|$ we get $\Hom_{K^{-}(\mc{N})}(X,N[i])= 0$. It follows immediately from the equivalence $\fc{\gamma}{K^b(\mc{T})}{D^b(\mc{U})}$ that for any $A,B\in D^b(\mc{U})$, $\Hom_{D^b(\mc{U})}(A,B[i])\cong 0$ for $i$ big enough. For any $i\in \mb{Z}$ and complex $X\in K^{-,b}(\mc{N})$ denote by $X^{\geq i}\in K^b(\mc{N})$ the complex $ X^{i} \rightarrow X^{i+1} \rightarrow \cdots \rightarrow X^{r_X} \rightarrow 0$. Note that for any $j>i$ we have,

$$ \Hom_{K^{-}(\mc{N})}(X,N[i]) \cong \Hom_{K^b(\mc{N})}(X^{\geq -j}, N[i]) \cong \Hom_{D^b(\mc{U})}(\gamma X^{\geq -j}, \gamma N[i]) $$
\\
Recall that $r_E$ has the property that for any $j<r_E$ the complex $X$ is $j$-acyclic (and so exact in this degree). Take $i<\min\{0,r_E\}$. Then the triangle $X^{\geq i-2} \rightarrow X \rightarrow \cok(d_X^{i-4})[-i+3] \rightarrow_{+1}$ and the fact that there are no extensions between tilting modules in the category $\mc{U}$ shows that $\Hom_{K^{-}(\mc{N})}(X, N[i]) \cong \Hom_{K^{-}(\mc{N})}(X^{\geq i-4}, N[i])$. But the later is zero for $i$ big enough.  So, $a_N$ is well-defined as well. The last statement now follows because $\Hom$ is a cohomological functor. \findem

The functions $a_N$ and $b_N$ for negligible tilting modules $N\in \mc{N}$ induce integer-valued functions on the Grothendieck group $K_0(K^{-,b}(\mc{N}))$. We denote this function by the same symbols, that is, we have functions $\fc{a_N,b_N}{K_0(K^{-,b}(\mc{N}))}{\mb{Z}}$ defined by $[X]\mapsto a_N(X)$ and $[X]\mapsto b_N(X)$, where $[X]$ is the image of $X$ in the Grothendieck group.\\

Assume $\mf{g}\neq \mf{sl}_2$. By Proposition \ref{closedretracts}, we have an exact sequence of groups

$$ \xymatrix{  K_0(K^b(\mc{N})) \ar[r] &  K_0(K^{-,b}(\mc{N})) \ar[r]^(0.4)p &  K_0(K^{-,b}(\mc{N})/K^b(\mc{N})) \ar[r] & 0}  $$

Let $V\in K^{-,b}(\mc{N})$ be an acyclic complex of the form

$$ \xymatrix{\cdots \ar[r] & N^{-3} \ar[r] & N^{-2} \ar[r] & P \ar[r] & N^0} $$

\noindent where $N^0$ is a negligible non projective, $P$ is its projective cover and $\cdots \rightarrow N^{-3} \rightarrow N^{-2}$ is an $\mc{N}$-resolution of $\ker(P\rightarrow N)^0$, it is infinite by Lemma \ref{globalNdim}.

\begin{prop} $p([V])\neq 0$. \end{prop}

\dem We have $b_{N^0}([V])\geq 1$ because at least we have the chain map given by the identity $N^0\rightarrow N^0$ in degree zero. We claim that $a_N([V])=0$ for any $N\in \mc{N}$. Indeed, using the notation in the proof of Lemma \ref{welldefiab} we have $ \Hom_{K^{-}(\mc{N})}(V, N[i]) =0$. If $p([V])=0$ we would have $[V]\in K_0(K^b(\mc{N}))$. Then, $[V]=\sum_{i=1}^m n_i[N_i]$ for some $N_i\in \mc{N}$ and $n_i\in \mb{Z}$. But since tilting modules are self dual this would give $a_N([V])=b_N([V])$ which is a contradiction when $N=N^0$.  \findem

Note that the functor $\fc{\overline{F}}{K^{-,b}(\mc{N})/K^b(\mc{N})}{D^b(\mc{U})/\langle \mc{N} \rangle}$ induces a surjective ring homomorphism $K_0(\overline{F}):  K_0(K^{-,b}(\mc{N})/K^b(\mc{N} )) \to K_0(D^b(\mc{U})/\langle \mc{N} \rangle)$.

\begin{prop}\label{enhancementofR} The ring homomorphism $K_0(\overline{F})$ is non-injective unless $\mf{g}=\mf{sl}_2$. \end{prop}

\dem If $\g\neq \mf{sl}_2$, there are negligible modules which are not projective, so we can construct a complex $V$ as above. By construction, $V\in K^{-,b,ex}(\mc{N})=\ker \overline{F}$. The induced map on Grothendieck groups $K_0(\overline{F})$ contains $p([V])$ in its kernel, and so $\ker(K_0(\overline{F}))\neq 0$. Hence, the Grothendieck group $K_0(K^{-,b}(\mc{N})/K^b(\mc{N}))$ is different and surjects to $K_0(D^b(\mc{U})/\langle \mc{N} \rangle)$. In the case of $\mf{g}=\mf{sl}_2$, $\mc{N}=\mc{P}$ and the functor $\overline{F}$ is an equivalence, then the map $K_0(\overline{F})$ is clearly injective. \findem

\begin{cor} If $\g\neq \mf{sl}_2$ then $K_0(K^{-,b,ex}(\mc{N}))\neq 0$. \end{cor}
\findem

We obtain a ring surjection $K_0({\bf S}(\mc{U}/\mc{N})) \to \mc{R}$ which is non-injective unless $\mf{g}=\mf{sl}_2$. Hence, $K_0({\bf S}(\mc{U})/\mc{N}))$ can be thought of as an enhancement of $\mc{R}$.

\section{Example: the case of $\mf{sl}_2$}

We analyze the Verdier quotient $K^b(\mc{T})/K^b(\mc{N})$ for the special case of $\mf{g}=\mf{sl}_2$. In this section,  $\mc{U}$ denotes the category of finite dimensional $\U_q$-modules of type 1 for the Lie algebra $\mf{sl}_2$. Denote by $\mc{P}$ and by $\mc{I}$ the categories of projective and injective objects in $\mc{U}$. Hence, $\mc{P}=\mc{I}$.\\

Let $m\in C_{\ell} = \{ 0,1,\ldots, \ell-2  \}$. For any $i\geq -1$, consider the weights: $m_{2i+1} = m + 2(i+1)\ell$ and $m_{2i}= 2\ell(i+1) - (m+2)$.\\

We have the following exact sequences for Weyl and tilting modules,

\begin{equation}\label{eq2} \xymatrix{ 0\ar[r] & \Delta(m_i) \ar[r] & T(m_i) \ar[r] & \Delta(m_{i-1}) \ar[r] & 0 } \end{equation}
\begin{equation}\label{eq3} \xymatrix{ 0\ar[r] & L(m_{i-1}) \ar[r] & \Delta(m_i) \ar[r] & L(m_i) \ar[r] & 0 } \end{equation}

\begin{lem}\label{projressl2} In $\mc{U}$, given a simple module $L(m)$, for $m\in C_{\ell}$, the following is a minimal projective resolution of it (i.e., an $\mc{N}$-resolution):
$$\xymatrix{ \cdots \ar[r] & P_{m_2} \ar[r] & P_{m_1} \ar[r] & P_{m_0} \ar[r] & L(m) \ar[r] & 0 }$$ where $P_{m_i}=T(m_i)$. \end{lem}

\dem Easily follows from the exact sequences \ref{eq2} and \ref{eq3}. \findem

Since the category $\mc{U}$ is a Frobenious category, $\mc{U}/\mc{P}$ is a triangulated category canonically equivalent to ${\bf S}(\mc{U}/\mc{P})$, see section \ref{contrfinitesubcat} fo the definition of the stabilization of a left triangulated category.

\begin{lem}\label{stablesl2} $K^b(\mc{T})/K^b(\mc{N}) \cong D^b(\mc{U})/\langle\mc{P}\rangle \cong {\bf S}(\mc{U}/\mc{P}) \equiv \mc{U}/\mc{P} = \U_q-\underline{\md}$.\end{lem}

\dem The negligible tilting modules coincides with the projective modules. Then by Theorem \ref{equivalenceingeneral} and section \ref{contrfinitesubcat} we obtain the result. \findem

\begin{prop}\label{lemmanonsemisimple} $K^b(\mc{T})/K^b(\mc{N})$ is generated as a triangulated category by the simple modules $L(m)$, $m\in C_{\ell}$. For any $m,n\in C_{\ell}$, $$\Hom_{K^b(\mc{T})/K^b(\mc{N})}(L(m),L(n)[k]) = \begin{cases} \mb{C} \mbox{  if   } m=n, k=0,-1 \\ 0 \mbox{  otherwise  } \end{cases}$$ \end{prop}

\dem The category is generated by $L(m)$, for $m\in C_{\ell}$, since they are the fusion modules. By lemma \ref{projressl2}, we have projective resolutions for $L(m)$ and $L(n)$, say $P_{\bullet}$ and $Q_{\bullet}$ respectively. By definition of the hom-spaces in the quotient category we have:

$$ \Hom_{K^b(\mc{T})/K^b(\mc{N})}(L(m),L(n)[k]) \cong {\underset{i,i-k\geq0}{\colim}} \Hom_{\mc{U}}(\Omega^{i}(L(m)),\Omega^{i-k}(L(n))) $$

\noindent where $\Omega^{j}(L(m))$ is the kernel of the map $P_{m_j}\rightarrow P_{m_{j-1}}$ in the projective resolution for $L(m)$. Similarly, for $L(n)$. By the exact sequence $(\ref{eq2})$ and by the basis of morphisms between tilting modules given in \cite{AST}, the morphisms $P_{m_j}\rightarrow P_{m_{j-1}}$ is a lifting of the surjection $P_{m_j} \twoheadrightarrow \Delta(m_{j-1})$. Thus, $\Omega^{j}(L(m)) \cong \Delta(m_j)$. If $m\neq n$, there are no morphisms between $\Delta(m_i)$ and $\Delta(n_{i-k})$ by the linkage principle. It remains to study the case when $m=n$.\\

Let $m=n$ and $k>0$. Then there is no vector of weight $m_i$ in $\Delta(m_{i-k})$. Hence, we will get that $\Hom_{\mc{U}}(\Delta(m_i),\Delta(m_{i-k})) = 0$. On the other hand, if $k\leq 0$ the weight $m_{i-k}$ is bigger than the weight $m_i$, and there are morphisms unique up to scalars $\Delta(m_i)\rightarrow \Delta(m_{i-k})$ just when $k=0$ or $k=-1$. When $k=0$ is the identity map and when $k=-1$ is the unique map with image $L(m_i)$. \\

By the exact sequences $(\ref{eq2})$ and  $(\ref{eq3})$, we see that the unique non-zero morphisms $\Delta(m_i)\rightarrow \Delta(m_i)$ and $\Delta(m_i)\twoheadrightarrow L(m_i) \hookrightarrow \Delta(m_{i+1})$, which are basis for $\Hom_{\mc{U}}(\Delta(m_i),\Delta(m_{i-k}))$, $k=0,-1$, does not factor trough a projective object. Thus, $\Hom_{\mc{U}}(\Delta(m_i),\Delta(m_{i-k}))\cong \mb{C}$. Moreover, the canonical maps $\Hom_{\mc{U}}(\Delta(m_i),\Delta(m_{i+1}))\rightarrow \Hom_{\mc{U}}(\Delta(m_{i+1}),\Delta(m_{i+2}))$ are isomorphisms. Therefore, $\Hom_{D^b(\mc{U}(\mf{sl}_2))/\langle \mc{N} \rangle}(L(m),L(m)[k]) \cong \mb{C}$ when $k=0$ or $k=-1$.
\findem

\section{The fusion ring of a spherical category and the small quantum
group}\label{sphericalsection}

In this section we recall the notion of a spherical category and suggest a definition of its fusion ring. As an example we partially describe the fusion ring of the small quantum group and shows that in the case of $\mf{sl}_2$ it coincides with the version of the Verlinde algebra introduced by Lachowska in \cite{AL}. We also discuss the problem of whether the restriction of negligible tilting modules is a contravariantly finite subcategory.

\subsection{Fusion rings for spherical categories}\label{sphercatdeffus}

Let $\mathscr{C}$ be a rigid monoidal category with unit object $\mb{1}$. We assume that $\mathscr{C}$ is a $\mb{k}$-linear category where $\mb{k}$ denotes the commutative ring $\End{\mb{1}}$. The category $\mathscr{C}$ is a pivotal category if it is endowed with a pivotal structure, that means, a monoidal isomorphisms between $X$ and $X^{**}$ for any object $X$ in $\mathscr{C}$. The pivotal structure implies that the right and left dualities coincide.\\

In a pivotal category $\mathscr{C}$ there are left and right traces $\fc{\Tr_L, \Tr_R}{\End(X)}{\mb{k}}$ for any $X\in \mathscr{C}$, see \cite{BW} for definitions.  For any two morphisms $f,g$ in $\mathscr{C}$ we have $\Tr_L(f\otimes g)=\Tr_L(f)\Tr_L(g)$. We say that the category $\mathscr{C}$ is spherical if it is a pivotal category in which the left and right traces coincide. In this case we define the categorical or quantum dimension of an object $X$ by $\dim_q(X)=\Tr_L(1_X)$.\\

Given a spherical category $\mathscr{C}$, it is possible to construct a quotient category of $\mathscr{C}$ which is spherical and semisimple, see Theorem 2.9 in \cite{BW}. Here we present some examples.

\begin{ex} The category $\mc{T}$ of tilting modules for a quantized enveloping algebra at a root of unity is an additive spherical category. For the subcategory of negligible tilting modules $\mc{N}$, the quotient $\mc{T}/\mc{N}=:\mc{F}$ is a spherical category with finitely many simples objects indexed by the weights in the principal alcove. \end{ex}
\begin{ex} The category $\mc{U}$ is a spherical abelian category. The quotient $\mc{U}/\mc{N}_{\mc{U}}$ is spherical too, but typically it has infinitely many isomorphism classes of simple objects. \end{ex}
\begin{ex} The category of representations for the small quantum group $\uu_q$ is a spherical category. This example is studied in detail in the next section. \end{ex}

\begin{defi} An object $X\in \mathscr{C}$ is called negligible if $\Tr_L(f)=0$ for any $f\in \End_{\mathscr{C}}(X)$. In particular, $\dim_q(X)=0$. \end{defi}

Let $\mathscr{C}$ be an abelian spherical category and let $\mathscr{N}_\mathscr{C}$ be its full subcategory of negligible objects. It is possible to prove that $\mathscr{N}_\mathscr{C}$ a tensor ideal. Consider the Verdier quotient $D^b(\mathscr{C})/\langle\mathscr{N}_{\mathscr{C}}\rangle$, where $\langle \mathscr{N}_\mathscr{C} \rangle$ is the full triangulated subcategory of $D^b(\mathscr{C})$ generated by the objects of $\mathscr{N}_\mathscr{C}$ and closed under retracts and tensor products with arbitrary modules.

\begin{defi}\label{fusionsphericaldefi} The derived fusion category of $\mathscr{C}$ is $D^b(\mathscr{C})/\langle\mathscr{N}_\mathscr{C}\rangle$ and its fusion ring $\mc{R}_{\mathscr{C}}$ is $K_0(D^b(\mathscr{C})/\langle\mathscr{N}_\mathscr{C}\rangle)$. \end{defi}

An advantage of this definition is that we do not need to define tilting modules in the spherical category in order to define the fusion ring. We just need the spherical structure for the definition of the category of negligible modules. In \cite{AAGITV} the still unsolved problem of how to define tilting modules for a spherical category is discussed.\\

For the rest of this section suppose $\mathscr{C}$ is an abelian Frobenious category with enough projective and injective objects. It is known that the projective objects belongs to $\mathscr{N}_\mathscr{C}$. We finish this section giving a characterization of the property that the retracts of the category $\langle\mathscr{N}_{\mathscr{C}}\rangle$ are generated by the objects of $\mathscr{N}_{\mathscr{C}}$.

\begin{lem}\label{conesstable} Let $N\in \mathscr{N}_\mathscr{C}$ and $\iota: N \rightarrow I_N$ be an injective hull of $N$. Then, $C_N:=\cok(\iota)\in\mathscr{N}_\mathscr{C}$. \end{lem}

\dem Consider the exact sequence  $ \xymatrix{0\ar[r]&N\ar[r]^{\iota}&I_N\ar[r]^{\pi}&C_N\ar[r]&0} $ and let $f:C_N\rightarrow C_N$ be any automorphisms of $C_N$. Because $\mathscr{C}$ is a Frobenious category, $I_N$ is projective and then there exists map $\tilde{f}:I_N\rightarrow I_N$ such that $f\pi=\pi\tilde{f}$. Let $\tilde{f}|_N$ the restriction of $\tilde{f}$ to $N$ and note that $\im(\tilde{f}|_N)\subset N$ because $\pi(\iota\tilde{f}|_N)=0$. Also, becuase $\tilde{f}$ is an extension of $f$ we have that $\im(\tilde{f}|_{C_N})\subset C_N$. Then, $\Tr_q(\tilde{f})=\Tr_q(\tilde{f})|_N + \Tr_q(f)$. It follows that $\Tr_q(f)=0$ becuase $N$ and $I_N$ are negligibles. So, $C_N\in \mathscr{N}_\mathscr{C}$.

\findem

By section \ref{contrfinitesubcat}, if $\mc{P}$ is the full subcategory of projective objects of $\mathscr{C}$ we have equivalences of categories ${\bf S}(\mathscr{C}/\mc{P}) \equiv \mathscr{C}/\mc{P} \equiv D^b(\mathscr{C})/ \langle \mc{P} \rangle$. Recall, that in the stable category $\mathscr{C}/\mc{P}$ the cone of a morphisms is constructed as follows: given a map $f:X\rightarrow Y$ in $\mathscr{C}/\mc{P}$ there is an injective map $(f,\iota_X): X \rightarrow Y\oplus I_X$ in $\mathscr{C}$, where $\iota_X: X\rightarrow I_X$ is the injective hull of $X$ and such that the class of $(f, \iota_X)$ is $f$. We define $\cone(f)$ as the pushout object given by the maps $f$ and $\iota_X$. Distinguished triangles in $\mathscr{C}/\mc{P}$ are the ones isomorphic to triangles of the form $X\rightarrow Y\rightarrow \cone(X\rightarrow Y)\rightarrow \cok(X\rightarrow Y)$.\\

Denote by $\underline{\mathscr{N}_\mathscr{C}}$ the image of $\mathscr{N}_\mathscr{C}$ in the stable category $\mathscr{C}/\mc{P}$ and by $\underline{Trg(\mathscr{N}_\mathscr{C})}$ the image of $Trg(\mathscr{N}_\mathscr{C})$ in the category $D^b(\mathscr{C})/ \langle \mc{P} \rangle$.

\begin{prop}\label{closedretractsequiv} $Trg(\mathscr{N}_\mathscr{C})$ is closed under retracts if and only if for any monomorphisms $f:N\rightarrow N'$ in $\mathscr{N}_\mathscr{C}$ we have $\cok(f)\in \mathscr{N}_\mathscr{C}$. \end{prop}

\dem Suppose $f:N\rightarrow N'$ in $\mathscr{N}_\mathscr{C}$ such that $\cok(f)\in \mathscr{N}_\mathscr{C}$. This condition and the Lemma \ref{conesstable} implies that $\underline{\mathscr{N}_\mathscr{C}}$ is closed under shifts and cones. So, under the equivalence $\mathscr{C}/\mc{P} \equiv D^b(\mathscr{C})/ \langle \mc{P} \rangle$ corresponds to $\underline{Trg(\mathscr{N}_\mathscr{C})}$. Now, assume we have an object $N$ of $\mathscr{C}/\mc{P}$ which is a direct summand of some $N'$ in $\underline{\mathscr{N}_\mathscr{C}}$, that means, $N$ is a direct summand of $N'\oplus P$ in the category $\mathscr{C}$ for some projective $P$. Because, $N'\oplus P$ is negligible we get that $N$ is too, and so, $N\in \underline{\mathscr{N}_\mathscr{C}}$. This implies that $\underline{Trg(\mathscr{N}_\mathscr{C})}$ is closed under retracts. Finally, because $\langle \mc{P} \rangle$ is closed under retracts we get that $Trg(\mathscr{N}_\mathscr{C})$ is closed under retracts in $D^b(\mathscr{C})$. The other direction in the proposition is clear. \findem

\subsection{Fusion for the small quantum group}\label{smallfusion}

The small quantum group $\uu_q$ is defined to be the subalgebra of $\U_q$ generated by $E_i, F_i$ and $K_i$.
It is a finite dimensional Hopf subalgebra of dimension $\ell^{\dim_{\mb{C}}\mf{g}}$. We denote the category of all integrable type 1 $\uu_q$-modules by $\mf{u}^{int}$ and by $\mf{u}$ its subcategory of finite dimensional modules.\\

It is known that $\mf{u}$ is a spherical category with the usual quantum trace. The category of negligible modules in $\mf{u}$ is denoted by $\mc{N}_{\mf{u}}$. The fusion category of the small quantum group is $D^b(\mf{u})/\langle \mc{N}_{\mf{u}} \rangle$ and its fusion ring
is $\mc{R}_{\mf{u}} = K_0(D^b(\mf{u})/\langle \mc{N}_{\mf{u}} \rangle)$. For $M \in \mc{U}$ denote by $M |_{\mf{u}}$
its restriction to $\mf{u}$. The simple objects of $\mf{u}$ are $L(\lambda) |_{\mf{u}}$, $\lambda \in P_\ell$ where $P_\ell = \{ \lambda\in P^{+} | \langle \lambda, \alpha^{\vee} \rangle <\ell,  \alpha \in \Delta \}$, see \cite{APW2}. From this it follows
that $\mf{u}$ (and $D^b(\mf{u})$ as a triangulated category) is generated by $\Delta(\lambda)|_{\mf{u}}$,  $\lambda \in P_\ell$ and
that $K_0(\mf{u}) = K_0(D^b(\mf{u}))$ is a free $\mb{Z}$-module with basis $[\Delta(\lambda)|_{\mf{u}}]$, $\lambda \in P_\ell$.
The restriction map $\mc{U} \to \mf{u}$ defines a surjective ring homomorphism $K_0(\mc{U}) \to K_0(\mf{u})$.\\

In \cite{AL}, Lachowska defined an algebra $\overline{Vr}:= \mb{C}\otimes_{\mb{Z}}\mc{R} \otimes_{K_0(\mc{U})} K_0(\mf{u})$ which is a counterpart for the small quantum group of the fusion ring (or in her terminology, Verlinde algebra) $\mb{C}\otimes_{\mb{Z}}\mc{R}$ of $\mc{U}$. Its representation theoretical meaning remains mysterious. She shows that $\overline{Vr} \cong \mb{C} \otimes_{\mb{Z}} K_0(\mf{u})/I$ where $I$ is the ideal generated by $[\Delta(\lambda)|_{\mf{u}}] + [\Delta(s_\alpha \bullet \lambda)|_{\mf{u}}]$, $s$ is a reflection in $W$, $\lambda \in P_\ell$ and that $\dim_{\mb{C}} (\overline{Vr})^{-} = (\dim_{\mb{C}} \mb{C}\otimes_{\mb{Z}}\mc{R})/(|P|/|Q|) = |\mc{X}|$. Here the $\bullet$-action is defined by $w\bullet\lambda =w\cdot\lambda\quad mod(\ell P)$, for $w\in W_{\ell}$, $\lambda \in P_\ell$, and $\mc{X}$ is the set of regular weights inside the fundamental domain for the $\bullet$-action, $ \overline{\mc{X}}$. A basis for $\overline{Vr}$ is $\{[L(\lambda) |_{\mf{u}}] = [\Delta(\lambda) |_{\mf{u}}]\}$, $\lambda \in \mc{X}$.

\begin{prop}\label{closedretractssmallsl2} Let $\mf{g}=\mf{sl}_2$, then $Trg(\mc{N}_{\mf{u}})$ is closed under retracts.\end{prop}

This proposition implies that $Trg(\mc{N}_{\mf{u}}) = \langle \mc{N}_{\mf{u}} \rangle$, and so all the retracts in $\langle \mc{N}_{\mf{u}} \rangle$ are negligible modules. We sketch the proof of this proposition in the next section.\\

\begin{prop}\label{lachcomp}\begin{enumerate} \item $\dim_\mb{C}  \mb{C}\otimes_{\mb{Z}}\mc{R}_{\mf{u}} \geq |\mc{X}|$. 
\item When $\g = \mathfrak{sl}_2$ we have a canonical ring isomorphism
$\mb{C}\otimes_{\mb{Z}}\mc{R}_{\mf{u}}\cong \overline{Vr}$.\end{enumerate}
\end{prop}

\dem (1) By the linkage principle for the small quantum group, \cite{AL}, proposition 2.7, we get
$D^b(\mf{u})/\langle \mc{N}_{\mf{u}} \rangle = \oplus_{\lambda \in
\overline{\mc{X}}} D^b(\mf{u})^\lambda/\langle
\mc{N}^\lambda_{\mf{u}} \rangle$ so that $$\mb{C}\otimes_{\mb{Z}}\mc{R}_{\mf{u}} =
\oplus_{\lambda \in \overline{\mc{X}}} {\mb{C}}\otimes_{\mb{Z}}K_0( D^b(\mf{u})^\lambda/\langle \mc{{N}}^\lambda_{\mf{u}}
\rangle).$$ Now for each $\lambda \in \mc{X}$ we have the
$\mb{C}$-linear map $\dim_q: {\mb{C}}\otimes_{\mb{Z}}K_0( D^b(\mf{u})^\lambda/\langle
\mc{{N}}^\lambda_{\mf{u}} \rangle) \to \mb{C}$ which is non-zero
since $\dim_q(\Delta(\lambda)|_{\mf{u}}) =
\dim_q(\Delta(\lambda))\neq 0$. This proves (1).

\medskip \noindent In the $\mathfrak{sl}_2$-case we have that
$\overline{Vr}\cong {\mb{C}}\otimes_{\mb{Z}}K_0(\mf{u})/I$ where $I$ is
the ideal generated by $[\Delta(\lambda)|_{\mf{u}}] + [\Delta(s
\bullet \lambda)|_{\mf{u}}]$, $\lambda \in P_\ell$, $s$ a reflection in $W$. But if $\lambda < s \bullet \lambda$ then
$\dim_q (\Delta(\lambda)|_{\mf{u}}) = - \dim_q (\Delta(s \bullet
\lambda)|_{\mf{u}})$ and there is an extension $E_\lambda \in
\operatorname{Ext}^1_{\mf{u}}(\Delta(\lambda)|_{\mf{u}},\Delta(s
\bullet \lambda)|_{\mf{u}})$ which can be described as follows: As
a $\mb{C}[K]$-module $E_\lambda = \Delta(\lambda)|_{\mf{u}} \oplus
\Delta(s \bullet \lambda)|_{\mf{u}}$; the action of $E$ and $F$ is
the same as it would be in the direct sum with the sole exception
that $F$ applied to the lowest weight vector of
$\Delta(s\bullet \lambda)|_{\mf{u}}$ equals the highest weight vector of
$\Delta(\lambda)|_{\mf{u}}$. Then $\dim_q E_\lambda = 0$
and since $E_\lambda$ is indecomposable we get $E_\lambda \in
\mc{N}_{\mf{u}}$. Thus we get a quotient map $D^b(\mf{u})/\langle
E_\lambda, \lambda \in \mc{X} \rangle \to D^b(\mf{u})/\langle
\mc{N}_{\mf{u}} \rangle$ and hence, (since $[E_\lambda] =
[\Delta(\lambda)|_{\mf{u}}] + [\Delta(s \bullet
\lambda)|_{\mf{u}}]$ and $\langle
\mc{N}_{\mf{u}} \rangle$ is closed under retracts by Proposition \ref{closedretractssmallsl2}) a surjective ring homomorphism
$$\overline{Vr} \cong {\mb{C}}\otimes_{\mb{Z}}K_0(D^b(\mf{u})/\langle
E_\lambda, \lambda \in \mc{X} \rangle) \to \mb{C}\otimes_{\mb{Z}}\mc{R}_{\mf{u}}.$$ This
is an isomorphism by (1). This proves (2).
\findem

\subsection{Sketch of the proof of proposition \ref{closedretractssmallsl2}}

The proof is a case by case analysis using the classification of the indecomposable modules for the small quantum group of $\mf{sl}_2$ given in \cite{CPr}. It would take several pages to write all the details and they are not enlightening. Therefore we opted for a brief sketch. The indecomposable modules in the category $\mf{u}$ come in five groups:

\begin{itemize}
\item Standard modules $\Delta(n)$ which are indecomposable for any $n$ not divisible by $\ell$.
\item Duals of standard modules, i.e., costandard modules, $\nabla(n)$.
\item Some maximal submodules of standard modules $C_{\lambda,\mu}(n) \subseteq \Delta(n)$. They can be understood in the following way: Consider an integer $n$ such that $\Delta(n)$ is indecomposable. Write $n=n_0 + \ell n_1$, where $0\leq n_0\leq \ell-1$. As a vector space $\Delta(n)$ is generated by vectors $v_0, \dots, v_{n_0}, v_{n_0+1}, \ldots,v_{n-n_0-1}, v_{n-n_0}, \ldots, v_n$ where they are ordered from the highest to the lowest weight vectors, i.e., $v_0$ and $v_n$ are the highest and lowest weight vectors respectively. Define $C_{\lambda,\mu}(n)$ to be generated by $\lambda v_0 + \mu v_{n-n_0}, \ldots, \lambda v_{n_0} + \mu v_n,  v_{n_0+1}, \ldots, v_{n-n_0-1}$, where $\lambda, \mu\in \mb{C}$ not both zero at the same time.
\item Duals of the maximal modules $C_{\lambda,\mu}(n)$.
\item Projective modules $P(n)$. For each $n\in P_{\ell}$ there is a simple module, its projective cover is $P(n)$.
\end{itemize}

The negligible modules are: projective, standard, costandard modules and the maximal submodules of standard modules which have $\mb{C}$-dimensions divisible by $\ell$. By Proposition \ref{closedretractsequiv} it is enough to prove that for any monomorphism $N\hookrightarrow N'$ of negligible modules its cokernel $N'/N$ is a negligible module.\\

Write $N=\oplus N_i$ where the $N_i$ are indecomposable negligible modules in the above list. Then by direct inspection we see that $N_i$ has a submodule $K$ of dimension $\ell$ so that $N_i/K$ is negligible. Thus $N'/N\cong (N'/K)/(N/K)$ so if we can prove that $N'/K$, $N/K$ are negligible we get by induction on $\dim_{\mb{C}}(N)$ that $N'/N$ is negligible as well.\\

Therefore we come to the case when $N$ is indecomposable and $\dim_{\mb{C}}(N)=\ell$. Moreover, if the map $f=(f_1, \ldots, f_n) : N \hookrightarrow N_1\oplus \cdots \oplus N_n$, an inspection process will show that one of the components of $f$ is injective, say $f_1$. A case by case analysis allow us to conclude that $\cok f$ is negligible.

\subsection{About the contravariantly finiteness of $\mc{N}_{\mf{u}}$}
We have adjoint pair of functors, restriction and induction,  $\xymatrix{\Res: \mc{U}^{int} \ar@<0.5ex>[r] &  \mf{u}^{int}:  \Ind \ar@<0.5ex>[l]}$, where $\Res(M)=M|_{\mf{u}}$ is the restriction functor  and $\Ind(V)=(\mc{O}_q(G)\otimes V)^{\uu_q}$ is the induction functor, see \cite{APW2}. Restriction is an exact functor and in this case the induction is also exact by theorem 4.8 in \cite{APW2}, so they induce an adjoint pair of functors on the level of derived categories, that we denote by the same symbols, $\Res:D^b(\mc{U}^{int}) \rightleftharpoons D^b(\mf{u}^{int}):\Ind$.\\

By the results in \cite{AG}, the functor $\Ind$ factors as follows. Let $(\mc{U}^{int}, \mc{O}(G))$ be the category of $\U_q$-equivariant $\mc{O}(G)$-modules. For $V\in \mf{u}^{int}$ let $\tilde{\Ind}(V)=\Ind(V)$ equipped with its natural $\mc{O}(G)$-module structure coming from the isomorphism $\mc{O}(G)\cong \mc{O}_q(G)^{\uu_q}$. Then $\fc{\tilde{\Ind}}{\mf{u}^{int}}{(\mc{U}^{int},\mc{O}(G))}$ becomes an equivalence of categories. In this situation $\Ind = for \circ \tilde{\Ind}$ where $\fc{for}{(\mc{U}^{int},\mc{O}(G))}{\mf{u}^{int}}$ is the functor that forgets the $\mc{O}(G)$-module structure.\\

Let $\langle D^b(\mf{u}^{int})_{sing} \rangle$ the smallest triangulated subcategory of $D^b(\mf{u}^{int})$ which contains $L(\lambda)|_{\mf{u}}$ for $\lambda\in P^{+}$ singular and is closed under retracts and tensor products with arbitrary modules.

\begin{lem} $\Ind(\langle D^b(\mf{u}^{int})_{sing} \rangle)\subseteq \langle D^b(\mc{U}^{int})_{sing} \rangle$ and $\Res(\langle D^b(\mc{U}^{int})_{sing} \rangle) \subseteq \langle D^b(\mf{u}^{int})_{sing} \rangle$. \end{lem}

\dem The second assertion is obvious. For the first assertion, note that $\langle D^b(\mf{u}^{int})_{sing} \rangle$ is generated by $L(\lambda)|_{\mf{u}}$, for $\lambda\in P^{+}$ singular, under triangles, shifts and tensor products with arbitrary modules. Therefore, it is enough to observe that $\Ind(L(\lambda)|_{\mf{u}}) = \mc{O}(G)\otimes L(\lambda) \in \langle D^b(\mc{U}^{int})_{sing} \rangle$. \findem

Denote by $\mc{N}|_{\mf{u}}$ the restriction of the category $\mc{N}$ to $\mf{u}$. We have the following partial result.

\begin{thm} If $V\in \mc{U}$ then $V|_{\mf{u}}$ has an $\mc{N}|_{\mf{u}}$-hull and an $\mc{N}|_{\mf{u}}$-approximation, which is the restriction of an $\mc{N}$-null and an $\mc{N}$-approximation of $V$. In particular, simple objects in the category $\mf{u}^{int}$ has $\mc{N}|_{\mf{u}}$-approximations. \end{thm}

\dem I prove it for $\N$-hulls, for $\N$-approximations is dual. Let $\mc{K}=K^{+}(\mc{N})$. Since $\mc{N}$ is dually $\mc{U}$-approximating we there is a $t$-structure on $\mc{K}$ with $\mc{K}^{\geq 0}$ consisting of complexes living in degrees $\geq 0$ and such that the truncation $\fc{\tau^{\geq 0}}{\mc{K}}{\mc{K}^{\geq 0}}$ is defined as the $\N$-hull of the naive truncation functor in $\mc{K}$ (see Proposition 3.10 of  \cite{AB}). It is left adjoint to the inclusion $\mc{K}^{\geq 0}\hookrightarrow \mc{K}$. Let $\mc{K}_{res} = K^{+}(\mc{N}|_{\mf{u}})$ and let $\mc{K}^{\geq 0}_{res} = \mc{K}_{res}\cap \mc{K}^{\geq 0}$.\\

Pick a finitely generated projective presentation $P^{-1}\rightarrow P^{0} \rightarrow V \rightarrow 0$ of $V$. Then the complex $[P^{-1}\rightarrow P^0]\in \mc{K}$. We define $$  \tau^{\geq 0}_{res}(\Res[P^{-1}\rightarrow P^0]) := \Res\tau^{\geq 0}[P^{-1}\rightarrow P^0]  $$
Then for $W\in \mc{K}^{\geq 0}_{res}$ we get

\begin{equation}\label{eq4}\begin{split} & \Hom_{\mc{K}^{\geq 0}_{res}}(\tau^{\geq 0}_{res}(\Res[P^{-1}\rightarrow P^0]), W) =  \Hom_{\mc{K}^{\geq 0}_{res}}(\Res\tau^{\geq 0}[P^{-1}\rightarrow P^0], W) \cong \\
& \Hom_{\mc{K}^{\geq 0}_{res}}(\tau^{\geq 0}[P^{-1}\rightarrow P^0], \Ind W) =  \Hom_{\mc{K}}([P^{-1}\rightarrow P^0], \Ind W)= \\
 & \Hom_{\mc{K}_{res}}(\Res[P^{-1}\rightarrow P^0], W) \end{split}\end{equation}

We write $\tau^{\geq 0}([P^{-1}\rightarrow P^0]) = [N^0\rightarrow N^1 \rightarrow \cdots]$. The cone of the natural transformation $Id\rightarrow \tau^{\geq 0}$ evaluated at $[P^{-1}\rightarrow P^0]$ gives the exact complex

$$ \xymatrix{P^{-1} \ar[r] & P^0 \ar[r] & N^0 \ar[r] & N^1 \ar[r] & \cdots } $$

In particular this gives an injective map $V=\cok(P^{-1}\rightarrow P^0)\rightarrow N^0$. We claim that the restricted map $V|_{\mf{u}}\rightarrow N^0|_{\mf{u}}$ is an $\mc{N}|_{\mf{u}}$-hull. To see this, take $W=N|_{\mf{u}}$ a restricted tilting module living in degree zero. Then we get from $(\ref{eq4})$ that

$$ \Hom_{\mf{u}_q}(V|_{\mf{u}}, N|_{\mf{u}}) = \Hom_{\mc{K}_{res}}(\Res[P^{-1}\rightarrow P^0], N|_{\mf{u}}) =  $$
$$ \Hom_{\mc{K}_{res}}([N^0|_{\mf{u}}\rightarrow N^1|_{\mf{u}} \rightarrow \cdots], N|_{\mf{u}}) $$

In particular, this shows that any map $V|_{\mf{u}}\rightarrow N|_{\mf{u}}$ factors trough our given map $V|_{\mf{u}})\rightarrow N^0|_{\mf{u}}$.

\findem

It would be interesting to know if the categories $\mc{N}|_{\mf{u}}$ or $\mc{N}_{\mf{u}}$ are contravariantly finite subcategories of the category $\mf{u}$, if the category $\mc{N}_{\mc{U}}$ is contravariantly finite subcategory of the category $\mc{U}$ and if the tensor ideals $\langle\mc{N}\rangle$ and $\langle\mc{N}_{\mc{U}}\rangle$ are the same.

\section{Acknowledgments} I am very grateful to my advisors Erik Backelin and Kobi Kremnizer for all the teachings and also for suggesting me to study this problem. I would like to thank Paul Bressler for useful conversations.

%\addcontentsline{toc}{chapter}{References}
  \bibliographystyle{plain}
  %\bibliographystyle{apalike}
  %\nocite{*}
 \bibliography{BiblioDFC}

\end{document}